\RequirePackage{amsmath} 
\documentclass{svjour3}
\usepackage[utf8]{inputenc}

\usepackage{amssymb}
\usepackage{mathtools}
\usepackage[normalem]{ulem}
\usepackage{hyperref}
\usepackage[usenames,dvipsnames]{color}
\hypersetup{colorlinks,
  citecolor=black,
 linkcolor=black,
 urlcolor=Blue}

\usepackage{float}
\usepackage{booktabs}
\usepackage{times}
\usepackage{graphicx}
\usepackage[final]{pdfpages}
\usepackage{longtable}
\usepackage{float}  
\usepackage{url}
\usepackage{xcolor}
\usepackage{tcolorbox}
\definecolor{pigpink}{HTML}{FDD7E4}
\definecolor{lcyan}{HTML}{E0FFFF}
\definecolor{mint}{HTML}{98FF98}

\newcommand{\eps}{\varepsilon}

\usepackage{marginnote}
\usepackage{soul}

\newcommand\fref[1]{Fig.~\ref{#1}}

\newcommand\Fref[1]{Figure~\ref{#1}}

\usepackage{multirow}
\usepackage{indentfirst}
\usepackage{todonotes}

\title{The uncoupling limit of identical Hopf bifurcations with an application to perceptual bistability}

\author{Alberto P\'erez-Cervera \and Peter Ashwin \and Gemma Huguet \and Tere M-Seara \and James Rankin}

\institute{ 
A. P\'erez-Cervera, G. Huguet and T. M-Seara \at
Departament de Matem\`atiques, Universitat Polit\`ecnica de
Catalunya, Avda. Diagonal 647, 08028 Barcelona (Spain), BGSMATH
\and
P. Ashwin, J. Rankin \at
Department of Mathematics, College of Engineering, Mathematics and Physical Sciences, University of Exeter, Harrison Building, North Park Rd, Exeter EX4 4QF, UK
 \\
EPSRC Centre for Predictive Modelling in Healthcare, University of Exeter, Exeter, EX4 4QJ, UK
}

\begin{document}



\maketitle

\begin{abstract}
We study the dynamics arising when two identical oscillators are coupled near a Hopf bifurcation where we assume a parameter $\epsilon$ uncouples the system at $\epsilon=0$. 
Using a normal form for $N=2$ identical systems undergoing Hopf bifurcation, we explore the dynamical properties. 
Matching the normal form coefficients to a coupled Wilson-Cowan oscillator network gives an understanding of different types of behaviour that arise in a model of 
perceptual bistability. Notably, we find bistability between in-phase and anti-phase solutions that demonstrates the feasibility for synchronisation to act as 
the mechanism by which periodic inputs can be segregated (rather than via strong inhibitory coupling, as in existing models).  
Using numerical continuation we confirm our theoretical analysis for small coupling strength and explore the bifurcation diagrams for large coupling strength, 
where the normal form approximation breaks down. 
\end{abstract}

\tableofcontents

\newpage
\textbf{List of abbreviations}\\

\begin{tabular}{l l}
	IP & In-phase\\
	AP & Anti-phase\\
	FP & Fixed Point\\
	LC & Limit Cycle\\
	PD & Period Doubling\\
	PF & Pitchfork Bifurcation\\
	TR & Torus Bifurcation\\
	LA & Low Amplitude\\
	HA & High Amplitude
\end{tabular}

\section{Introduction}\label{sec:intro}

The Hopf bifurcation is a generic and well-characterized transition that a nonlinear system can undergo to create temporal patterns of behaviour on changing a parameter. 
At such a bifurcation, an equilibrium of an autonomous smooth dynamical system develops an oscillatory instability and emits a small amplitude periodic orbit that, 
when followed, may be used to understand a wide variety of oscillatory phenomena. 
This includes many problems that appear in Neuroscience applications \cite{AshComNic2016}.

For larger network systems composed of similar subsystems that undergo oscillatory instability, 
when coupled together  this can lead to the formation of non-trivial spatio-temporal patterns. 
Notably there is a large literature on coupled oscillators, viewed from a wide variety of theoretical view points, 
and from the point of view of applications, e.g. \cite{Pikovsky01}. 
Much of this theory either considers very specific models, or makes an assumption of weak coupling which allows a reduction to a phase oscillator description 
such as that of Kuramoto \cite{Acebron2005}, suitable for answering a lot of questions about synchronization of system oscillations.

In this paper we consider identical subsystems undergoing a Hopf bifurcation that have an uncoupling limit. 
This approach gives a natural setting of two parameters that allows a thorough and generic analysis of the low-dimensional dynamics of coupled oscillator systems, 
by means of normal form theory. 
We use this analysis to understand the behaviour of a pair of Wilson-Cowan oscillators that arise in a model of perceptual bistability, which complements the results in \cite{Borisyuk1995}.

The phenomenon of perceptual bistability motivates this study of
oscillatory dynamics in a coupled dynamical system. For certain static
but ambiguous sensory stimuli, two distinct perceptual interpretations
(percepts) are possible, but only one can be held at a time. Not only
can the initial percept be different from one short presentation of
the stimulus to the next, but for extended presentations, the percept
can switch dynamically. Perceptual bistability has been investigated
in a number of different visual paradigms e.g. ambiguous
figures~\cite{necker:32,rubin:21}, binocular
rivalry~\cite{levelt1968binocular,blake:89,blake:01}, random-dot
rotating spheres~\cite{wallach-oconnell:53},
motion plaids~\cite{hupe2003dynamics} and multistable barber-pole
illusion~\cite{meso2016relative}. Such ambiguous stimuli provide an
opportunity to gain insights about the computations underlying
perceptual competition in the brain. Whilst synchrony of oscillatory activity is known to play a role in the encoding of perceptually ambiguous 
stimuli~\cite{fries1997synchronization}, this mechanism has been widely overlooked.

Further background and motivation for the study of coupled oscillatory instabilities close to the uncoupling limit is given in Section~\ref{sec:intromaths}, 
whilst further background and motivation for the study of oscillatory dynamics in the context of perceptual competition is given in Section~\ref{sec:introneuro} 
(not required reading if primarily interested in this paper's mathematical results).

\subsection{Coupled oscillatory instabilities}\label{sec:intromaths}

As noted by several authors, networks of oscillators near Hopf bifurcation allow one to explore not just the collective phase dynamics but also amplitude 
behaviour \cite{golubitsky2003symmetry} and this allows one to use many of the tools of generic bifurcation theory with symmetry 
(in particular, the consequence of group actions on normal forms and the phase space) to understand the creating and properties of many oscillator patterns that may arise, 
biological applications including, for example, animal gaits and visual hallucination patterns \cite{golubitsky2003symmetry}. 

A recent paper \cite{ashwin2016hopf} explored coupled Hopf bifurcations in a two parameter setting where one of the parameters results in uncoupling of the systems. 
In that setting, they found that it is possible to find not only a reduction to Kuramoto-like oscillators in a weak coupling close to threshold limit, 
but also to find the next order corrections that include multiple oscillator interactions. 
The setting also allows study of patterns where only part of the system is oscillating. 
More precisely, \cite{ashwin2016hopf} considers $N$ identical and identically interacting smooth ($C^{\infty}$) vector fields on $x_i \in \mathbb{R}^d$ $(d \geq 2)$ and 
present the normal form
near a Hopf bifurcation. 
  
In this paper we explore the dynamical properties of the special case $N=2$ with $d=2$. 
We give a dimension reduction via group-invariant coordinates in order to simplify dynamics.  
In the 2D normal form we look at the effects of coupling beyond the weak limit. A similar analysis was performed in \cite{Aronson1990} for the case of a linear coupling term, thus considering a particular sub-case of the normal form 
studied here.
We then apply this theory to understand the appearance of a variety of oscillatory patterns in a model of perceptual bistability. 

We emphasize that we explore a special case of two identical Hopf bifurcations that has symmetries and is close to $1:1$ resonance. The case of a double Hopf bifurcation without symmetries has been studied in \cite{gavrilov1980some} (see also \cite{guckenheimer2013nonlinear}); by assuming non-resonant conditions on the Hopf bifurcation frequencies, the author provides a normal form for the bifurcation and performs a detailed study of the dynamics. Depending on the value of the coefficients very rich dynamics can be found. However, our study examines different behaviours, that are generic for systems with symmetries close to identical Hopf bifurcations, but not in the more general case.

\subsection{Oscillatory models of perceptual bistability}\label{sec:introneuro}

\begin{figure}[t]
\centering
{\includegraphics[width=120mm]{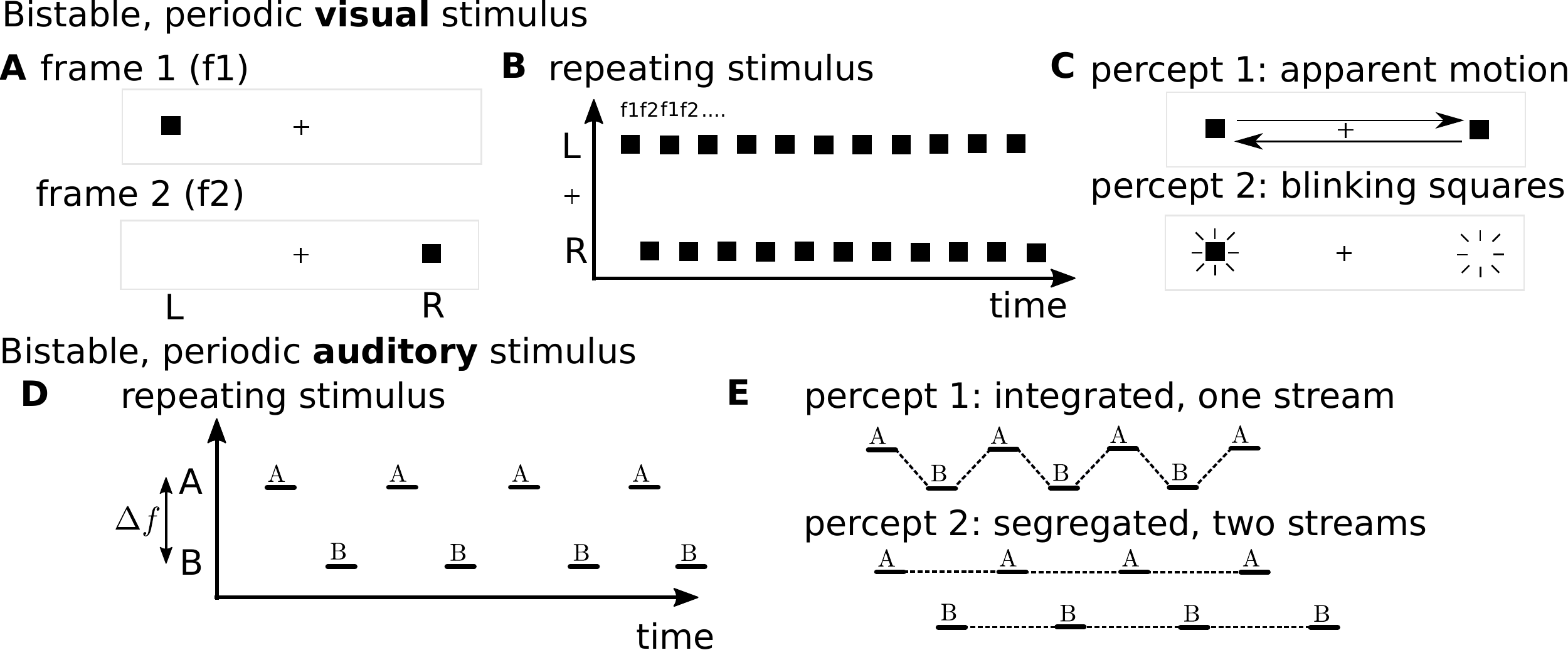}}
\caption{Perceptually bistable stimuli that repeat
  periodically. \textbf{A}: The visual stimulus alternates between two
  frames f1 and f2 with a black square flipping between a position to
  the left (L) or right (R) of a fixed fixation point
  ($+$). \textbf{B}: Schematic of repeating stimulus, illustrating
  that the black square alternates between the L (f1) and R (f2)
  locations over time. \textbf{C}: Two perceptual interpretations are
  possible for the same repeating stimulus: 1) apparent motion where a
  single square appears to travel between the L and R locations, 2)
  blinking squares, where the squares appear to separately blink on
  and off at two distinct locations. \textbf{D}: The auditory stimulus
  features alternating pure tones at frequencies A and B separated by
  a difference in frequency $\Delta f$. \textbf{E}: Two perceptual
  interpretations are possible for the same repeating stimulus: 1) an
  integrated percept where both the A and B tone sequences are heard
  in a single alternating stream ``ABABAB\ldots'', 2) a segregated
  percept where the A tone sequence ``A{\_}A{\_}A{\_}\ldots'' is heard
  separately from the B tone sequence ``{\_}B{\_}B{\_}B\ldots''
  (``{\_}'' is a silent gap).  } \label{fig:schematic}
\end{figure}

Perceptual bistability can also arise with stimuli that change
periodically. Apparent motion can be observed when a dot on a screen
present at one location disappears and spontaneously reappears at a
nearby location, as if travelling smoothly across the
screen~\cite{kolers1964illusion,anstis1980perception}. \Fref{fig:schematic}A
shows two frames of such an apparent motion display\footnote{More complex example than the one we're interested in:
  \href{https://open-mind.net/videomaterials/kohler-motion-quartet.mp4/view}{https://open-mind.net/videomaterials/kohler-motion-quartet.mp4/view}},
where a black square to the left of a fixation point might reappear on
the right of the fixation point. If two such frames alternate every,
say 200\,ms, as in the schematic \fref{fig:schematic}B, this can be
perceived as a single square moving from side to side (``percept 1''
in \fref{fig:schematic}C). However, another interpretation is
possible, of distinct squares blinking on and off either side of the
fixation point (``percept 2'' in \fref{fig:schematic}C). Watching such
a display, perception switches between percept 1 and percept 2 every
few seconds;
see~\cite{anstis1985adaptation-b,ramachandran1983perceptual},
references within and more
recently~\cite{gilroy2004multiplicative,muckli2002apparent}. Perceptual
bistability also occurs for the so-called auditory streaming
paradigm~\cite{van1975temporal,anstis1985adaptation,Pressnitzer2006}\footnote{\href{https://auditoryneuroscience.com/scene-analysis/streaming-alternating-tones}{https://auditoryneuroscience.com/scene-analysis/streaming-alternating-tones}}.
The stimulus consists of interleaved sequences of tones A and B,
separated by a difference in tone frequency $\Delta f$, and repeating
in an ``ABABAB\ldots'' pattern (\fref{fig:schematic}D). This can be
perceived as one stream, integrated into an alternating rhythm
(``percept 1'' in \fref{fig:schematic}E) or as two segregated streams
(``percept 2'' in \fref{fig:schematic}E); see recent reviews
\cite{moore2002factors,snyder2017recent}. There are commonalities
between these visual and auditory paradigms, in percept 1
(\fref{fig:schematic}C and E) the stimulus elements are linked into a
single percept. In percept 2, the stimulus elements are separated into
their distinct parts in space or in frequency. In both cases the
stimulus alternates rapidly (in the range at 2--5\,Hz for the visual
stimulus~\cite{anstis1985adaptation-b}; in the range 5--10\,Hz for the
auditory stimulus~\cite{van1975temporal}), whilst the perceptual
interpretations are stable on the order of several seconds (over many
cycles of the rapidly alternating stimuli).

Models of perceptual bistability have successfully captured the
dynamics of perceptual switching
\cite{Laing2002,Wilson2003,wilson2007minimal}, the dependence of these
dynamics on stimulus parameters
\cite{Laing2002,moreno2010alternation,seely2011role,rankin2015neuromechanistic},
mechanisms for attention \cite{li2017attention}, entrainment to slowly
varying stimuli \cite{kim2006stochastic} and the effects of stimulus
perturbations~\cite{rankin2017stimulus}. Generally models are based on
competition between abstract, percept-based units
\cite{Wilson2003,Shpiro2009,huguet2014noise,li2017attention}, but more
recently models with a feature-based representation of competition
have been
developed~\cite{Laing2002,kilpatrick2013short,rankin2014bifurcation,rankin2015neuromechanistic}. Some
percept-based models have explored how rapidly alternating inputs
($>$2\,Hz) can still give rise to stable perception over several
seconds
\cite{wilson2007minimal,vattikuti2016canonical,li2017attention}. The
models described above have considered competition directly between
populations encoding different percepts, or between populations
separated on a feature space. In general model studies of perceptual bistability have not explored how synchrony properties of oscillations entrained at the rate of a rapidly
alternating stimulus could be the mechanism by which different
perceptual interpretations emerge and coexist as bistable
states (although see \cite{wang2008oscillatory} for a large network approach to this problem). We hypothesise that oscillations play a key role in
perceptual integration (such as ``percept 1'') and perceptual segregation
(such as ``percept 2''). Towards exploring this hypothesis in future
modelling studies of perceptual bistability, this paper lays the
mathematical groundwork for studying the encoding perceptual states
similar to those described above. An aim of the study is to identify
regions of parameter space where such states coexist for a suitable
neural oscillator model (but not transitions between these states). 

Matching the normal form
  coefficients to a coupled Wilson-Cowan oscillator network allows for
  an understanding of the parameters in the model that govern different
  types of behaviour. Numerical continuation is used to confirm our
  theoretical analysis and to complete bifurcation diagrams for large
  coupling strength demonstrating where the normal form approximation
  breaks down. Finally, our analysis is extended with numerics to
  demonstrate that coexisting states akin to ``percept 1'' and ``percept
  2'' persist in the presence of symmetrical periodic inputs. These coexisting states persist with low coupling strengths (down to the uncoupling limit) thus removing the need for the assumption of strong mutual inhibition between neural populations encoding different perceptual interpretations.

\subsection{Outline}
The structure of the paper is as follows: in Section \ref{sec:section2} we use recent theoretical results in \cite{ashwin2016hopf} 
to write the normal form of a system of two weakly coupled identical oscillators near a Hopf bifurcation. 
In Section \ref{sec:section3} we perform a dynamical analysis of the system given by the dominant terms of the normal form. 
In particular, we study how the solutions for the uncoupled system persist for weak coupling. 
In Section \ref{sec:section4} we identify different dynamical regimes depending on specific coefficients of the normal form and study the bifurcation diagrams. 
In Section~\ref{sec:percept} we write the equations for two mutually inhibiting Wilson-Cowan oscillators near a Hopf bifurcation and we perform a change of 
coordinates to put the system in the normal form discussed in Section \ref{sec:section2}. 
For this example, we compare the theoretical predictions given by the normal form analysis with a bifurcation diagram computed numerically. 
Finally, we note that the results are of broad interest, extending beyond the study of neural oscillators and perceptual bistability to the study of any 
system involving two coupled oscillators.

\section{Two identical Hopf bifurcations with an uncoupling limit}\label{sec:section2}

We will study systems consisting of two identically coupled oscillators of the form:
\begin{equation}\label{eq:originalSystem}
\begin{gathered}
\frac{d x_1}{dt} = H_\lambda(x_1) + \epsilon h_{\lambda, \epsilon}(x_1; x_2), \\
\frac{d x_2}{dt} = H_\lambda(x_2) + \epsilon h_{\lambda, \epsilon}(x_2; x_1), 
\end{gathered}\quad \quad x_1, x_2 \in  \mathbb{R}^2 \quad \quad \epsilon, \lambda \in  \mathbb{R}
\end{equation}
having $S_2$ permutation symmetry. 
We assume that when system \eqref{eq:originalSystem} is uncoupled ($\epsilon=0$), each system undergoes a Hopf bifurcation at the origin when the parameter 
$\lambda$ crosses zero. 

More concretely, we assume that the uncoupled system for $x \in \mathbb{R}^2$ given by
\begin{equation*}
\frac{dx}{dt} = H_\lambda(x)
\end{equation*}
has a stable focus at $x=0$ for $\lambda < 0$ that undergoes a supercritical Hopf bifurcation for $\lambda = 0$
which gives rise to a small amplitude stable limit cycle for $\lambda > 0$.
For simplicity  we assume that the eigenvalues of $DH_{\lambda}(0)$ are $\lambda \pm i\omega$ with $\omega \neq  0$. Moreover, without loss of generality, we assume that  
$(x_1, x_2) = (0,0)$ is an equilibrium point for $(\lambda, \epsilon)$ in some neighbourhood of $(0,0)$ for system \eqref{eq:originalSystem}.

\subsection{Truncated Normal Form in Complex Coordinates}

In \cite{ashwin2016hopf}, it is shown that systems as in \eqref{eq:originalSystem}, having  $S_2$ symmetry and undergoing a supercritical Hopf 
bifurcation for $\lambda = 0$, can be written in the following normal form
%
\begin{equation}\label{eq:anotherEquation}
\begin{gathered}
\frac{dz_1}{dt} = U_\lambda(z_1) + \epsilon F_N(z_1, z_2, \epsilon) + \mathcal{O}_{N+1}(z_1, z_2),\\
\frac{dz_2}{dt} = U_\lambda(z_2) + \epsilon F_N(z_2, z_1, \epsilon) + \mathcal{O}_{N+1}(z_2, z_1),
\end{gathered}\quad \quad z_1, z_2 \in  \mathbb{C} 
\end{equation}
where $F_N$ is a $N$-degree polynomial function that is equivariant under the  
rotational symmetries 
$$
\big(F_N(z_1 e^{i\phi}, z_2 e^{i\phi}, \epsilon) = e^{i\phi} F_N(z_1, z_2, \epsilon)\big).
$$
If we consider the normal form up to order 
three and ignore the $\mathcal{O}_4(z)$ terms, we obtain the truncated normal form
\begin{equation}\label{eq:tereSystem}
\begin{split}
\frac{dz_1}{dt} &= z_1 \Big(\lambda + i\omega + \alpha_{01}|z_1|^2 \Big) + \epsilon \Big[ z_1\Big(\alpha_{\epsilon0} + \alpha_{\epsilon1}|z_1|^2 
+ \alpha_{\epsilon2}|z_2|^2 + \alpha_{\epsilon3}\bar{z}_2z_1\Big) \\ &+ z_2\Big(\beta_{\epsilon0} 
+ \beta_{\epsilon1}|z_1|^2 + \beta_{\epsilon2}|z_2|^2 + \beta_{\epsilon3}\bar{z}_1z_2\Big)  \Big] 
, \\
\frac{dz_2}{dt} &= z_2 \Big(\lambda + i\omega + \alpha_{01}|z_2|^2 \Big) + \epsilon \Big[ z_2\Big(\alpha_{\epsilon0} + \alpha_{\epsilon1}|z_2|^2 
+ \alpha_{\epsilon2}|z_1|^2 + \alpha_{\epsilon3}\bar{z}_1z_2\Big) \\ &+ z_1\Big(\beta_{\epsilon0} + \beta_{\epsilon1}|z_2|^2 + \beta_{\epsilon2}|z_1|^2 
+ \beta_{\epsilon3}\bar{z}_2z_1\Big)  \Big] 
,
\end{split}
\end{equation}
where the constants $\alpha_{01}, \alpha_{\epsilon i}, \beta_{\epsilon i} \in \mathbb{C}$ with the restriction  $Re(\alpha_{01}) < 0$ 
because the Hopf bifurcation is supercritical.

\section{Dynamical analysis of the truncated normal form}\label{sec:section3}

\subsection{Hopf bifurcations of the origin}\label{sec:fixPointSec}

It is straightforward to check that the origin
\begin{equation}\label{eq:fixedPointSolution}
\mathcal{S}_0 = \Big\{ z_1 = z_2 = 0 \Big\},
\end{equation}
is a fixed point of the normal form \eqref{eq:tereSystem}. Let us start by analysing its stability. The Jacobian matrix of system \eqref{eq:tereSystem} evaluated at the origin is
\begin{equation}\label{eq:jacobianMatrix}
\left(\begin{array}{cccc} \lambda + i\omega + \epsilon \alpha_{\epsilon0} & 0 & \epsilon \beta_{\epsilon0} & 0\\
0 & \lambda - i\omega + \epsilon \bar{\alpha}_{\epsilon0}  & 0 & \epsilon \bar{\beta}_{\epsilon0} \\
\epsilon \beta_{\epsilon0} & 0 & \lambda + i\omega + \epsilon \alpha_{\epsilon0} & 0 \\
0 & \epsilon \bar{\beta}_{\epsilon0} & 0 & \lambda - i\omega + \epsilon \bar{\alpha}_{\epsilon0}  \\
\end{array}\right),
\end{equation}
and their eigenvalues are given by 
\begin{equation}\label{eq:theEigenvalues}
\begin{array}{rcl}
\mu_+ = \lambda + i\omega + \epsilon(\alpha_{\epsilon0} + \beta_{\epsilon0}), \quad \quad
\mu_- = \lambda + i\omega + \epsilon(\alpha_{\epsilon0} - \beta_{\epsilon0}),  \end{array}
\end{equation}
and its complex conjugate pairs ($\bar \mu_+, \bar \mu_-$). 

Clearly, when $\epsilon = 0$ the origin undergoes a double Hopf bifurcation at $\lambda = 0$. More interestingly, for $\epsilon \neq 0$, the origin undergoes two independent Hopf bifurcations, 
given by $Re (\mu_+)=0$ and $Re (\mu_-)=0$. These conditions define the following Hopf bifurcation curves $\mathcal{C}^{\pm}_{HB}$ in the ($\lambda$, $\epsilon$)-parameter space
\begin{equation}\label{eq:hbCurves}
\begin{aligned}
C^+_{HB} &= \Big\{ Re(\mu_+) = 0 \quad \text{or equivalently} \quad \bar{\alpha}^+ := \lambda + \epsilon(\alpha_{\epsilon0R} + \beta_{\epsilon0R}) = 0 \Big\}, \\
C^-_{HB} &= \Big\{ Re(\mu_-) = 0 \quad \text{or equivalently} \quad \bar{\alpha}^- := \lambda + \epsilon(\alpha_{\epsilon0R} - \beta_{\epsilon0R}) = 0 \Big\}.
\end{aligned}
\end{equation}
At each curve $C^{\pm}_{HB}$, a limit cycle is born, that will be denoted by $\mathcal{S}^{\pm}_{osc}$.

To study the stability of the origin of system \eqref{eq:tereSystem}, we analyse the sign of the real part of its eigenvalues 
$\mu^+$ and $\mu^-$ given in \eqref{eq:theEigenvalues} at the Hopf bifurcation curves $C^{\pm}_{HB}$ defined in \eqref{eq:hbCurves}. Thus,
\begin{equation}
\begin{aligned}
\text{if} \quad (\lambda, \epsilon) \in C^+_{HB} \rightarrow & \quad Re(\mu_+) = 0, \quad Re(\mu_-) = -2 \epsilon \beta_{\epsilon0R}, \\
\text{if} \quad (\lambda, \epsilon) \in C^-_{HB} \rightarrow & \quad Re(\mu_+) = 2 \epsilon \beta_{\epsilon0R}, \quad Re(\mu_-) = 0 .
\end{aligned}
\end{equation}

\begin{figure}[H]
	\centering
	{\includegraphics[width=110mm]{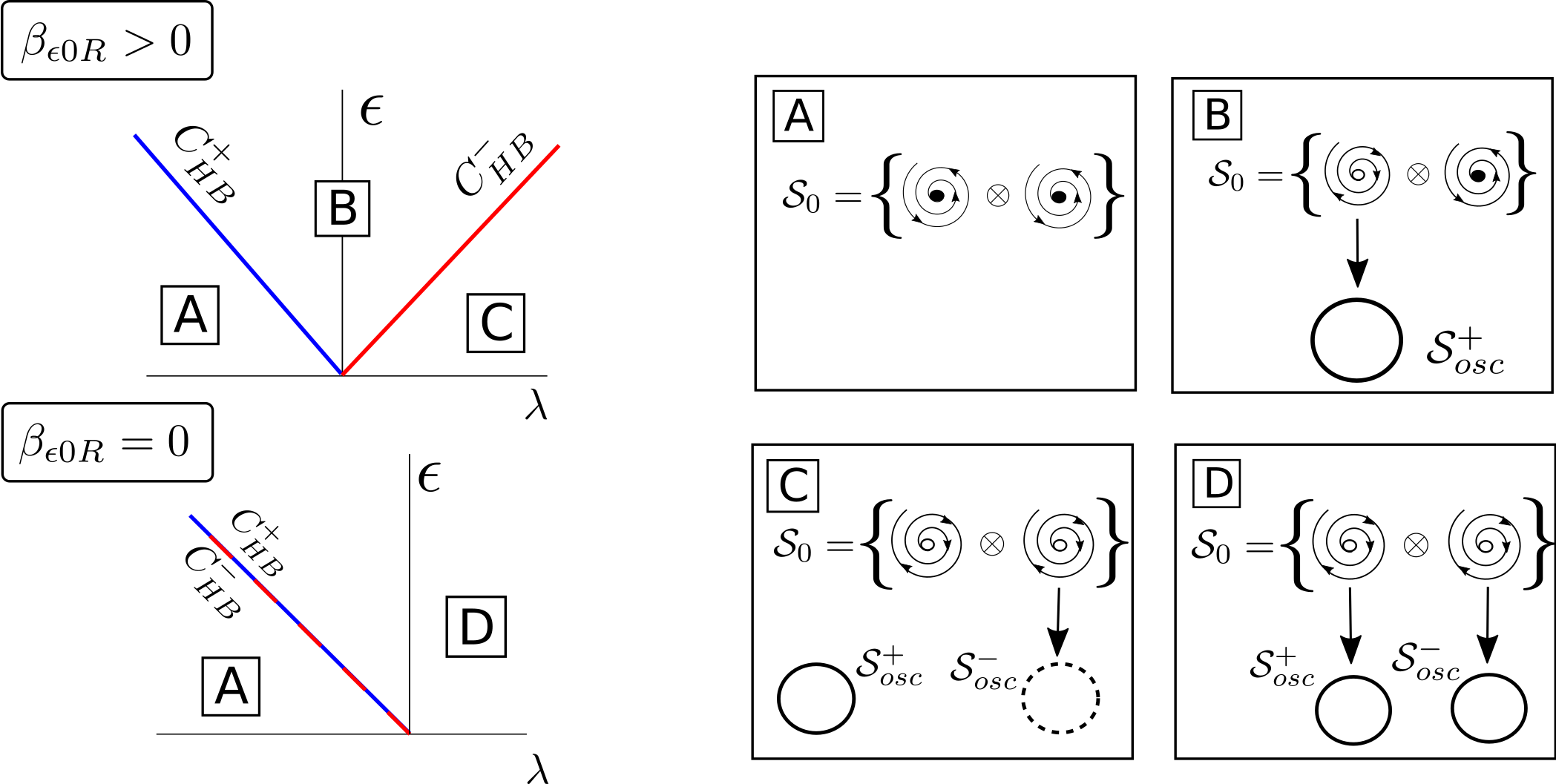}}
	\caption{Sketch for the curves $C^{\pm}_{HB}$ in \eqref{eq:anotherEquation}. 
		If $\beta_{\epsilon0R} > 0$ a stable limit cycle emerges from $C^{+}_{HB}$ whereas a saddle limit cycle emerges from $C^{-}_{HB}$. The case $\beta_{\epsilon0R} < 0$ is analogous just reversing $\pm$ by $\mp$. 
		For the special case $\beta_{\epsilon0R} = 0$, two stable limit cycles emerge at the coincident curves $C^{+}_{HB}$ and $C^{-}_{HB}$. For these plots we assume  $\beta_{\epsilon0R} > \alpha_{\epsilon0R} > 0$.} \label{fig:nova2}
\end{figure}

Therefore, we conclude that (see Fig. \ref{fig:nova2}):
\begin{itemize}
	\item 
	If $\beta_{\epsilon0R} > 0$, for $(\lambda, \epsilon) \in C^{+}_{HB}$ the solution $\mathcal{S}_0$ changes from a stable focus to a saddle-focus and a stable limit cycle $\mathcal{S}_{osc}^+$ emerges from $C^{+}_{HB}$. 
	Moreover, when $(\lambda, \epsilon) \in C^{-}_{HB}$, the solution $\mathcal{S}_0$ changes from a saddle-focus to an unstable focus and a saddle limit cycle 
	$\mathcal{S}_{osc}^-$ appears. 
	\item 
	If $\beta_{\epsilon0R} < 0$, for $(\lambda, \epsilon) \in C^{-}_{HB}$ the solution $\mathcal{S}_0$ changes from a stable focus to a saddle-focus and a stable limit cycle $\mathcal{S}_{osc}^-$ emerges from $C^{-}_{HB}$. 
	Moreover, when $(\lambda, \epsilon) \in C^{+}_{HB}$, the solution $\mathcal{S}_0$ changes from a saddle-focus to an unstable focus and a saddle limit cycle 
	$\mathcal{S}_{osc}^+$ appears. 
	\item If $\beta_{\epsilon0R} = 0$,   for $(\lambda, \epsilon) \in C^{-}_{HB}
	=C^{+}_{HB}$, the
	solution $\mathcal{S}_0$ changes from a stable focus to an unstable focus and two stable limit cycles 
	$\mathcal{S}_{osc}^+$ and $\mathcal{S}_{osc}^-$, appear. 
\end{itemize}

In the next section we analyse the oscillatory solutions $\mathcal{S}_{osc}^\pm$ that arise from the bifurcation curves $C^{\pm}_{HB}$ of system \eqref{eq:tereSystem}.


\subsection{Truncated Normal Form in Polar Coordinates}
\label{sec:polarCoords}

To perform the analysis of the oscillatory solutions $\mathcal{S}_{osc}^\pm$ we express the normal form in \eqref{eq:tereSystem} in polar coordinates, that is, we write 
$z_n = r_n e^{i\varphi_n}$ with $r_n > 0$ and $\varphi_n \in \mathbb{T}$:
\begin{equation}\label{eq:polarEqs}
\begin{split}
\dot{r}_1 &= r_1 \Big( \lambda + \alpha_{01R} r^2_1 \Big) + \epsilon f_r(r_1, r_2, \Delta \varphi),\\
\dot{r}_2 &= r_2 \Big( \lambda + \alpha_{01R} r^2_2 \Big) + \epsilon f_r(r_2, r_1, -\Delta \varphi),\\
r_1 \dot{\varphi}_1 &= r_1 \Big( \omega + \alpha_{01I} r^2_1 \Big) + \epsilon f_\varphi(r_1, r_2, \Delta \varphi),\\
r_2 \dot{\varphi}_2 &= r_2 \Big( \omega + \alpha_{01I} r^2_2 \Big) + \epsilon f_\varphi(r_2, r_1, -\Delta \varphi),
\end{split}
\end{equation}
where $\Delta \varphi = \varphi_2 - \varphi_1$ and the subscript $X = {R, I}$ in $\alpha_{01}$ refers to its real and imaginary parts, respectively. The expression for the functions $f_r$ and $f_\varphi$ can be found in Eq. \eqref{eq:basicEqs} in the Appendix.
%
System \eqref{eq:polarEqs} can be also written using the variable $\Delta \varphi$:
\begin{equation}\label{eq:phaseDifEqs}
\begin{aligned}
\dot{r}_1 &= r_1 \Big( \lambda + \alpha_{01R} r^2_1 \Big) + \epsilon f_r(r_1, r_2, \Delta \varphi),\\
\dot{r}_2 &= r_2 \Big( \lambda + \alpha_{01R} r^2_2 \Big) + \epsilon f_r(r_2, r_1, -\Delta \varphi),\\
\dot{\Delta \varphi} &= \alpha_{01I} (r^2_2 - r^2_1) + \epsilon f_{\Delta \varphi}(r_1, r_2, \Delta \varphi),\\
\dot{\varphi}_1 &= \omega + \alpha_{01I} r^2_1 + \frac{\epsilon}{r_1 } f_\varphi(r_1, r_2, \Delta \varphi),
\end{aligned}
\end{equation}
where the expression for the function $f_{\Delta \varphi}$ can be found in Eq. \eqref{eq:phaseDifTerms} in the Appendix.

\begin{remark}
The general non-resonant case of the double Hopf bifurcation is discussed in \cite{gavrilov1980some} (see also \cite{guckenheimer2013nonlinear}). The equations for the normal form in polar coordinates satisfy that the amplitudes $r_1, r_2$ decouple from the angles $\varphi_1, \varphi_2$. However, in our case (see system \eqref{eq:phaseDifEqs}) the equations for the amplitudes $r_1, r_2$ depend on $\Delta \varphi = \varphi_2 - \varphi_1$ leading to different generic dynamics than the one in \cite{gavrilov1980some}, which we study in this paper.
\end{remark}

Notice that the analysis of system \eqref{eq:phaseDifEqs} can be simplified by studying the system consisting of the first three equations, since they can be decoupled from the last one. Furthermore, we can further simplify the analysis by exploiting the $S_2$ permutation symmetry of the system. This symmetry acts on phase space as
\begin{equation}\label{eq:invariantMap}
K: \ (r_1, r_2, \Delta \varphi) \rightarrow (r_2, r_1, -\Delta \varphi) \quad \text{and} \quad K^2 = Id.
\end{equation}
This action can be diagonalised using sum and difference variables $s = r_1 + r_2$, $d = r_1 - r_2$, with $s, d \in \mathbb{R}^+\times\mathbb{R}$: in this case
\begin{equation}\label{eq:invariantEqiMap}
\tilde{K}(s, d, \Delta \varphi) \rightarrow (s, -d, -\Delta \varphi).
\end{equation}
Thus, expressing the first three equation of system \eqref{eq:phaseDifEqs} in the variables $(s, d, \Delta \varphi)$ we have
\begin{equation}\label{eq:sdEqs}
\begin{split}
\dot{s} &= s \Big(\lambda + \frac{\alpha_{01R}}{4}(s^2 + 3d^2)\Big) + \epsilon g_s(s, d, \Delta \varphi),\\
\dot{d} &= d \Big(\lambda + \frac{\alpha_{01R}}{4}(d^2 + 3s^2)\Big) + \epsilon g_d(s, d, \Delta \varphi),
\\
\dot{\Delta \varphi} &= -\alpha_{01I}sd + \epsilon g_{\Delta \varphi}(s, d, \Delta \varphi),
\end{split}
\end{equation}
where the expressions for functions $g_s$, $g_d$ and $g_{\Delta \varphi}$ are given in Eq. \eqref{eq:sdFunctions} in the Appendix.

The system \eqref{eq:sdEqs} will be the object of study for the rest of the Section and will be referred to as the \textit{reduced system}. As we will see in Section \ref{sec:oscSol}, working in the variables $s, d, \Delta \varphi$ has the advantage that the linearised system about the solutions of interest becomes block diagonal.

\subsubsection{Dynamical analysis of the reduced system in the uncoupled case ($\epsilon$ = 0)}\label{sec:uncSD}

The general picture of the uncoupled case can be obtained straightforwardly from the original system  \eqref{eq:tereSystem} for $\epsilon=0$. Indeed, as we consider two identical systems having a supercritical Hopf bifurcation at $\lambda = 0$, 
the solutions of system \eqref{eq:tereSystem} for $\lambda>0$ will correspond to the cartesian 
product of solutions of each $2$-dimensional system. In this Section we show how the solutions for $\epsilon = 0$ are seen in the reduced system \eqref{eq:sdEqs} so that we can explore how they evolve for $\epsilon \neq 0$. System \eqref{eq:sdEqs} for $\epsilon = 0$ writes
\begin{equation}\label{eq:uncoupledSD}
\begin{split}
\dot{s} &= s \Big(\lambda + \frac{\alpha_{01R}}{4}(s^2 + 3d^2)\Big),\\
\dot{d} &= d \Big(\lambda + \frac{\alpha_{01R}}{4}(d^2 + 3s^2)\Big),\\
\dot{\Delta \varphi} &= - \alpha_{01I}sd.
\end{split}
\end{equation}

Notice that in this case, the first two equations decouple from the third one and can be studied independently. As the variables $(s, d)$ are defined in 
$\mathbb{R}^+\times\mathbb{R}$, the fixed points of the first two equations of system \eqref{eq:uncoupledSD} are given by
\begin{equation}\label{eq:thePoints}
\begin{aligned}
\Big(0, 0 \Big), \quad 		\Big(\sqrt{\frac{-4\lambda}{\alpha_{01R}}}, 0 \Big), \quad 
\Big(+\sqrt{\frac{-\lambda}{\alpha_{01R}}},-\sqrt{\frac{-\lambda}{\alpha_{01R}}} \Big),  \quad 
\Big(+\sqrt{\frac{-\lambda}{\alpha_{01R}}},+\sqrt{\frac{-\lambda}{\alpha_{01R}}}\Big). 
\end{aligned} 
\end{equation}
Then, as the Jacobian matrix for the two first equations of system \eqref{eq:uncoupledSD} is given by
\begin{equation}\label{eq:matrixRedSys}
\left(\begin{array}{cc} \lambda + \frac{3\alpha_{01R}}{4}(s^2 + d^2) & \frac{\alpha_{01R}}{4}6ds\\ 
\frac{\alpha_{01R}}{4}6ds & \lambda + \frac{3\alpha_{01R}}{4}(s^2 + d^2)  \end{array}\right),
\end{equation}
it is straightforward to see that the eigenvalues of \eqref{eq:matrixRedSys} for 
$(s, d) = (0, 0)$
are $\lambda$ (double), for ($s, d$) = $\Big(\sqrt{\frac{-4\lambda}{\alpha_{01R}}}, 0\Big)$ 
are $-2\lambda$ (double) and for ($s, d$) = $\Big(\sqrt{\frac{-\lambda}{\alpha_{01R}}}, \pm \sqrt{\frac{-\lambda}{\alpha_{01R}}}\Big)$ 
are $\lambda$ and $-2 \lambda$.

Thus, when $\lambda = 0$ the origin undergoes a bifurcation and changes from stable to unstable while three new fixed points appear: one stable corresponding to 
$(s, d) = \Big( \sqrt{\frac{-\lambda}{\alpha_{01R}}}, 0\Big)$ plus two unstable corresponding to  
$(s, d) =  \Big( \sqrt{\frac{-\lambda}{\alpha_{01R}}}, \pm \sqrt{\frac{-\lambda}{\alpha_{01R}}}\Big)$.

\begin{figure}[H]
	\centering
	{\includegraphics[width=100mm]{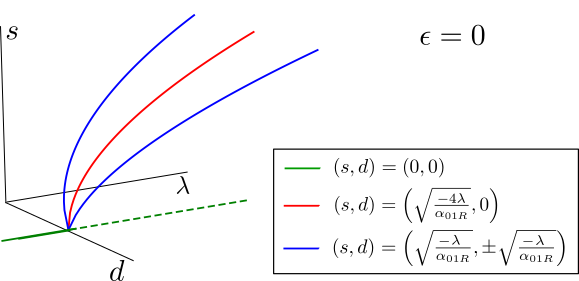}}
	\caption{Bifurcation diagram of system \eqref{eq:uncoupledSD} for $\epsilon=0$ as a function of $\lambda$. For the critical value $\lambda = 0$ the system undergoes a  bifurcation.}\label{fig:unperturbedSD1}
\end{figure}

Now let us study the solutions of system \eqref{eq:uncoupledSD} obtained from the fixed points \eqref{eq:thePoints} when considering the variable $\Delta \varphi$. The (singular) solution 
\begin{equation}
\bar {\mathcal{S}}_0 = \Big\{ s = d = 0, \quad \quad \Delta \varphi \in  \mathbb{T} \Big\},
\end{equation}
corresponds to the origin $\mathcal{S}_0$ in \eqref{eq:fixedPointSolution} of system \eqref{eq:tereSystem}, which is a focus with eigenvalues $\lambda \pm i\omega$ (double) (see Section \ref{sec:fixPointSec}).

For any value $\Delta\varphi_0$, the solution 
\begin{equation}
\bar {\mathcal{S}}_1 (\Delta\varphi_0) = \Big\{ s = \sqrt{\frac{-4\lambda}{\alpha_{01R}}}, \quad \quad d = 0, \quad \quad \Delta \varphi =\Delta \varphi_0 \Big\},
\end{equation}
is a fixed point of system \eqref{eq:uncoupledSD} with eigenvalues $-2\lambda$ (double) and $0$. These fixed points fill up the invariant curve
\begin{equation}\label{eq:theTorus3D}
\bar {\mathcal{T}}_0 = \Big\{ s = \sqrt{\frac{-4\lambda}{\alpha_{01R}}}, \quad \quad d = 0, \quad \quad \Delta \varphi \in \mathbb{T} \Big\},
\end{equation}
whose characteristic exponents are $-2\lambda$ (double). 
The fixed points $\bar {\mathcal{S}}_1 (\Delta\varphi_0)$ and the invariant curve $\bar {\mathcal{T}}_0$  correspond in the original system \eqref{eq:tereSystem} for $\epsilon = 0$ to the periodic orbits
\begin{equation}
\begin{aligned}
\mathcal{S}_1(\varphi_2^0) = 
\Big\{z_1=\sqrt{\frac{-\lambda}{\alpha_{01R}}}e^{i\varphi_1(t)}, \quad & \quad z_2=\sqrt{\frac{-\lambda}{\alpha_{01R}}}e^{i(\varphi_1(t) + \Delta \varphi_0)}, \\
 \Delta \varphi_0 = \varphi^0_2 - \varphi^0_1, \quad & \quad 
\varphi_1(t) = \varphi^0_1 + \big(\omega - \lambda \frac{\alpha_{01I}}{\alpha_{01R}}\big)t \quad,\quad t \in \mathbb{R} \Big\},
\end{aligned}
\end{equation}
and the $2$-dimensional invariant torus $\mathcal{T}_0$
\begin{equation}\label{eq:theTorus}
\mathcal{T}_0 = \bigcup_{\varphi_2^0 \in \mathbb{T}} \mathcal{S}_1(\varphi_2^0)=
\Big\{ |z_1|=|z_2|= \sqrt{\frac{-\lambda}{\alpha_{01R}}}, \quad \quad \varphi_1, \varphi_2 \in \mathbb{T} \Big\},
\end{equation}
respectively. Notice that the periodic orbits $\mathcal{S}_1(\varphi_2^0)$ fill the torus $\mathcal{T}_0$. The characteristic exponents of $\mathcal{T}_0$ are the eigenvalues of the fixed point $(s, d) = (\sqrt{\frac{-4\lambda}{\alpha_{01R}}}, 0)$ of the first two equations of system \eqref{eq:uncoupledSD} which are $-2\lambda$ (double). 

The invariant 2-torus  $\mathcal{T}_0$ is the product of two periodic orbits with the same period in the uncoupled case $\epsilon = 0$. Note that $\mathcal{T}_0$ is normally hyperbolic as each periodic orbit is linearly stable and the torus is foliated with periodic orbits; see for example \cite{ashwin2016hopf}. We recall that roughly speaking an invariant manifold is normally hyperbolic if the dynamics in the normal directions expands or contracts at a stronger rate than the internal dynamics. In our case the normal dynamics near the torus behaves as $e^{-2\lambda t}$ whereas the internal dynamics is just a rotation. Therefore the torus $\mathcal{T}_0$ is normally hyperbolic.

The last two fixed points in \eqref{eq:thePoints} give rise to the following periodic orbits of system \eqref{eq:uncoupledSD}
\begin{equation}
\begin{aligned}	
\bar {\mathcal{S}}^2 = 
\Big\{ s = d= \sqrt{\frac{-\lambda}{\alpha_{01R}}}, \quad \quad & \quad \Delta \varphi= \Delta \varphi_0- \frac{\alpha_{01I}}{\alpha_{01R}}\lambda t, \quad t \in \mathbb{R} 
\Big\},\\
\bar {\mathcal{S}}^3 = \Big\{ s =- d= \sqrt{\frac{-\lambda}{\alpha_{01R}}},  \quad & \quad \Delta \varphi= \Delta \varphi_0+ \frac{\alpha_{01I}}{\alpha_{01R}}\lambda t, \quad t \in \mathbb{R}  \Big\},
\end{aligned}
\end{equation}
whose characteristic exponents are $\lambda$ and $-2\lambda$, so they are of saddle type. 
These solutions correspond to the periodic solutions 
\begin{equation}
\begin{array}{rcl}	
\mathcal{S}^2 &=& \Big\{ z_1= \sqrt{\frac{-\lambda}{\alpha_{01R}}}e^{i\varphi_1(t)}, 
\quad z_2 = 0, \quad \varphi_1(t) = \varphi^0_1 + \big(\omega - \lambda \frac{\alpha_{01I}}{\alpha_{01R}}\big)t, \quad t \in \mathbb{R} \Big\},\\
\mathcal{S}^3 &=& \Big\{ z_1= 0, \quad z_2 = \sqrt{\frac{-\lambda}{\alpha_{01R}}}e^{i\varphi_2(t)}, \quad \varphi_2(t) 
= \varphi^0_2 + \big(\omega - \lambda \frac{\alpha_{01I}}{\alpha_{01R}}\big)t, \quad t \in \mathbb{R} \Big\},
\end{array}
\end{equation}
of the original system \eqref{eq:tereSystem} for $\epsilon = 0$ which have characteristic exponents $-2\lambda, \lambda\pm i\omega$. 
Therefore, they are hyperbolic periodic orbits of saddle type for $\lambda > 0$.

\begin{figure}[H]
	\centering
	{\includegraphics[width=60mm]{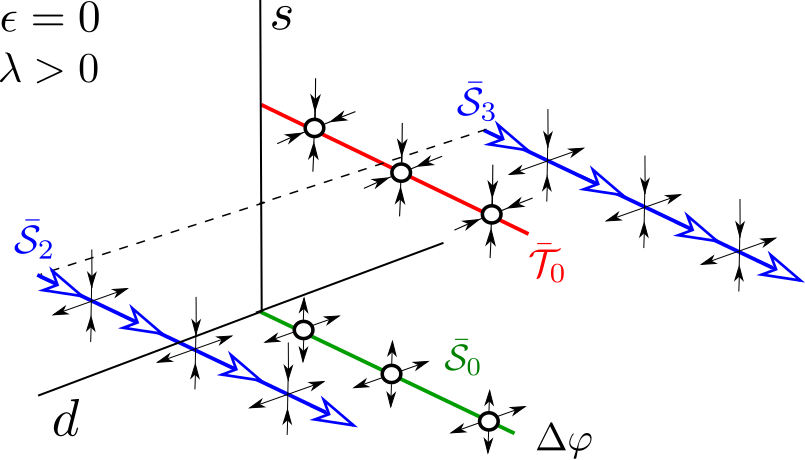}}
	\caption{Phase space for the unperturbed system \eqref{eq:uncoupledSD} for $\lambda > 0$. There are two invariant curves, $\bar {\mathcal{S}}_0$ (which is unstable) and $\bar {\mathcal{T}}_0$ (which is stable), filled with fixed points. 
		Moreover there exist two saddle periodic orbits $\bar {\mathcal{S}}_2$ and $\bar {\mathcal{S}}_3$.}\label{fig:unperturbedSD2}
\end{figure}

In conclusion, for $\epsilon = 0$, the 4D solutions $\mathcal{S}_0$, $\mathcal{T}_0$ and $\mathcal{S}^2$ and $\mathcal{S}^3$ arising from the union of solutions of each independent subsystem in \eqref{eq:tereSystem} can be seen in the uncoupled reduced system \eqref{eq:uncoupledSD} as two invariant curves filled with fixed points, $\bar {\mathcal{S}}_0$ and $\bar {\mathcal{T}}_0$,
and two saddle periodic orbits, $\bar{\mathcal{S}}^2$ and $\bar{\mathcal{S}}^3$, respectively (see Fig. \ref{fig:unperturbedSD2}).

Solutions $\mathcal{S}^2$ and $\mathcal{S}^3$ are hyperbolic periodic orbits for $\lambda > 0$ and  $\epsilon=0$.
Therefore,  for $\lambda > 0$ fixed and $\epsilon$ small enough there exist periodic orbits $\mathcal{S}^{2}_{\epsilon }$ and $\mathcal{S}^{3}_{\epsilon}$ that are  $C^1$-close to the unperturbed ones.

To ensure the persistence of the torus  $\mathcal{T}_0$ we use the Fenichel theorem \cite{Fenichel71} which guarantees the persistence of normally hyperbolic invariant manifolds (with a certain degree of smoothness) for small enough perturbations.
%
\begin{lemma}
	For a fixed value of $\lambda>0$, there exists  $\epsilon_0$ = $\epsilon_0(\lambda)$, such that for any $0\le \epsilon\le \epsilon_0$, 
	system \eqref{eq:tereSystem} has a stable $2$-dimensional torus  $\mathcal{T}_\epsilon$ 
	that is  $\mathcal{C}^1$-close to $\mathcal{T}_0$.
	\label{lem:persistence}
\end{lemma}

The analytic continuation when $\epsilon$ increases of the periodic orbits $\mathcal{S}_\epsilon^2$, $\mathcal{S}_\epsilon^3$ 
and the invariant torus $\mathcal{T}_{\epsilon}$ provided by Lemma~\ref{lem:persistence} is beyond the scope of this paper.
We note that the periodic orbits $\mathcal{S}_\epsilon^2$ and $\mathcal{S}_\epsilon^3$ are limited only by hyperbolicity.
Moreover, previous work \cite{ashwin2016hopf} highlighted that continuation of the torus $\mathcal{T}_{\epsilon}$ with $\epsilon$  in Lemma~\ref{lem:persistence} is only possible for $\epsilon=o(\lambda)$. Beyond this regime there will typically be loss of smoothness and breakup of the torus \cite{afraimovich1991invariant}.

In Section~\ref{sec:oscSol} we are able to study the persistence, for $(\lambda,\epsilon)$ small, of the periodic solutions $\mathcal{S}_{osc}^{\pm}$ that are born at the bifurcation curves $C_{HB}^{\pm}$ (see~\eqref{eq:hbCurves}). In Remark~\ref{rem:rem2} we relate these periodic orbits $\mathcal{S}_{osc}^{\pm}$ with the invariant torus $\mathcal{T}_{\epsilon}$
for $\lambda$ fixed and $\epsilon$ small enough, i.e. where the existence of the invariant torus is guaranteed. Later, in Section~\ref{sec:section4} we give a detailed study of all the
possible bifurcations of the solutions $\mathcal{S}_{osc}^{\pm}$.

\subsubsection{The oscillating solutions $\mathcal{S}_{osc}^{\pm}$ in the coupled case ($\epsilon > 0$)}\label{sec:oscSol}


We can take advantage of the $S_2$ symmetry of system \eqref{eq:sdEqs} to look for solutions which remain invariant under 
the application of the permutation map $\tilde {K}$ in \eqref{eq:invariantEqiMap}. 
Notice that by denoting $r_1 = r_2 = r^*$, the curves $(r^*, r^*, 0)$ and $(r^*, r^*, \pi)$ are invariant for system \eqref{eq:sdEqs}. 
Then, if we write these curves in the ($s,d$) coordinates 
\begin{equation}
\Xi^+ = \Big\{ (s, d, \Delta \varphi) = (s, 0, 0) \Big\}, \quad \quad \quad
\Xi^- = \Big\{ (s, d, \Delta \varphi) = (s, 0, \pi) \Big\},
\end{equation}
the dynamics for system \eqref{eq:sdEqs} when restricted to $\Xi^\pm$ reduces to 
\begin{equation}\label{eq:equivariantSub}
\begin{split}
\dot{s} =& \lambda s + \frac{s^3\alpha_{01R}}{4} + \epsilon \Big[(\alpha_{\epsilon0R} \pm \beta_{\epsilon0R})s + \\
&
\frac{s^3}{4}(\overbrace{\alpha_{\epsilon2R} + \alpha_{\epsilon1R} + \beta_{\epsilon3R} \pm (\beta_{\epsilon2R} + \beta_{\epsilon1R} +
\alpha_{\epsilon3R})}^{K_{stb}^\pm} ) \Big],\\
\dot{d} =& 0, \\
\dot{\Delta \varphi} =& 0,
\end{split}
\end{equation} 
where the $\pm$ sign corresponds to $\Delta \varphi = 0, \pi$, respectively.

It is straightforward to check that the equation for $s$ in \eqref{eq:equivariantSub} has three steady solutions, namely, 
$s=0$ (which corresponds to the solution $\bar{\mathcal{S}}_0$ studied before) and  $s^{\pm}_{osc}$ given by
\begin{equation}\label{eq:equivariantSolutions}
s^{\pm}_{osc} = \sqrt{ \frac{-4\left(\lambda+ \epsilon (\alpha_{\epsilon 0R} \pm  \beta_{\epsilon0R}) \right)}{\alpha_{01R} + \epsilon K_{stb}^\pm}}.
\end{equation} 
Notice that since $s \in \mathbb{R}^+$ we have discarded the negative solutions for the square root.

Taking into account that $\alpha_{01R}<0$, solutions $s^\pm_{osc}$ in \eqref{eq:equivariantSolutions} are only admissible when
$\bar{\alpha}^\pm = \lambda + \epsilon(\alpha_{\epsilon0R} \pm \beta_{\epsilon0R}) > 0$. This restriction defines the following conditions for the bifurcation
\begin{equation}\label{eq:lesAlfes}
\begin{split}
\bar{\alpha}^+ = \epsilon(\alpha_{\epsilon0R} + \beta_{\epsilon0R}) + \lambda = 0 & \quad \quad \text{for} \quad \Delta \varphi = 0,\\
\bar{\alpha}^- = \epsilon(\alpha_{\epsilon0R} - \beta_{\epsilon0R}) + \lambda = 0 & \quad \quad \text{for} \quad \Delta \varphi = \pi,
\end{split} 
\end{equation}
which are exactly the conditions defining the curves $C^{\pm}_{HB}$ in \eqref{eq:hbCurves} corresponding to the Hopf bifurcations of the origin. 

Therefore, for ($\lambda, \epsilon$)-values on the right-hand-side of curves $C^{\pm}_{HB}$ we can define, respectively, 
the following fixed points of system \eqref{eq:sdEqs}
\begin{equation}\label{eq:theSoscSols}
\begin{split}
\bar{\mathcal{S}}^+_{osc} = (s, d, \Delta \varphi) = (s^+_{osc}, 0, 0), \\
\bar{\mathcal{S}}^-_{osc} = (s, d, \Delta \varphi) = (s^-_{osc}, 0, \pi),
\end{split}
\end{equation}
which appear across a pitchfork bifurcation (whose character will be discussed below) of the origin in the $s$ direction. Fixed points in \eqref{eq:theSoscSols} correspond to the periodic orbits $\mathcal{S}^{\pm}_{osc}$ of system \eqref{eq:tereSystem} that appear at the Hopf bifurcation curves. Next, we will study its 
stability and possible bifurcations by using the reduced system \eqref{eq:sdEqs}. 

The Jacobian matrix  evaluated at the fixed points $\bar{\mathcal{S}}^\pm_{osc}$ is block diagonal 
\begin{equation}\label{eq:jacobianMatrixSD}
\left(\begin{array}{ccc} c^s_s & 0 & 0\\ 0 & c^d_d & c^d_{\Delta\varphi}\\ 0 & c^{\Delta\varphi}_d & c^{\Delta\varphi}_{\Delta\varphi} \end{array}\right), 
\end{equation}
where the terms $c^s_s$, $c^d_d$,  $c^d_{\Delta\varphi}$, $c^{\Delta\varphi}_d$ and  $c^{\Delta\varphi}_{\Delta\varphi}$ are different from zero, and their precise expressions 
are given in Eq. \eqref{eq:jacobianTermsCs} in the Appendix.

\begin{figure}[H]
	\centering
	{\includegraphics[width=100mm]{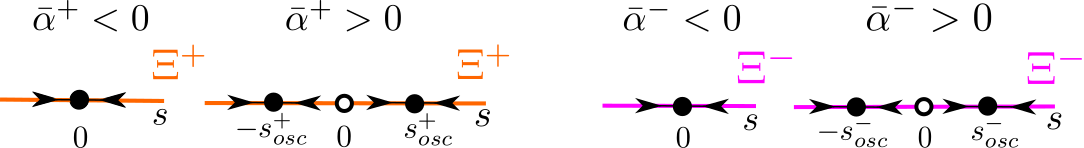}}
	\caption{Solutions $s^\pm_{osc}$ appear through a supercritical pitchfork bifurcation of the origin in the $s$ direction which takes place at the critical value $\bar{\alpha}^\pm = 0$ of the bifurcation parameter $\bar{\alpha}^\pm = \lambda + \epsilon(\alpha_{\epsilon 0R} \pm \beta_{\epsilon 0R}$).}\label{fig:perturbedSDpitch}
\end{figure}

Because of the block diagonal form of the Jacobian matrix, it is straightforward to check the stability in the $s$ direction as it corresponds to the 1x1 block. Thus,
the eigenvalue $\bar{\mu}^{\pm}_1$ takes the form
\begin{equation}\label{eq:sEigen}
\bar{\mu}^{\pm}_1 = c_s^s=-2(\epsilon(\alpha_{\epsilon0R} \pm \beta_{\epsilon0R}) + \lambda),
\end{equation}
and therefore, the solutions $\bar{\mathcal{S}}^\pm_{osc}$ are always stable in the $s$ direction as they appear for 
$\bar{\alpha}^\pm = \epsilon(\alpha_{\epsilon0R} \pm \beta_{\epsilon0R}) + \lambda > 0$. Therefore, the pitchfork bifurcations of the origin are supercritical (see Fig. \ref{fig:perturbedSDpitch}).

As the solutions $\bar{\mathcal{S}}^\pm_{osc}$ are always stable in the $s$ direction, one has to consider the eigenvalues of the 2x2 block, corresponding to the transverse directions, in order to study possible bifurcations of the symmetric solutions $\bar{\mathcal{S}}^\pm_{osc}$. 
The trace ($Tr^{\pm}$) and the determinant ($Det^{\pm}$) of the 2x2 block of 
\eqref{eq:jacobianMatrixSD} at $\bar{\mathcal{S}}^\pm_{osc}$ are given up to order 2 in $\lambda, \epsilon$ by: 
\begin{eqnarray}\label{eq:traceTereMatrix}
Tr^{\pm}(\lambda, \epsilon) &=& c^d_d + c^{\Delta\varphi}_{\Delta\varphi} = -2 \big(\lambda + \epsilon(\alpha_{\epsilon0R} \pm 3\beta_{\epsilon0R}) \big),\\
\label{eq:theDeterminant}
Det^\pm(\lambda, \epsilon) &=& \pm 4\epsilon \big(\lambda + \epsilon(\alpha_{\epsilon0R} \pm \beta_{\epsilon0R}) \big)(C_{det} + \beta_{\epsilon0R}) + 4\epsilon^2(\beta^2_{\epsilon0I} + \beta^2_{\epsilon0R}),
\end{eqnarray}
where 
\begin{equation}
C_{det} := \frac{\beta_{\epsilon0I}\alpha_{01I}}{\alpha_{01R}}.
\end{equation}

So, computing the discriminant
\begin{equation}
\label{eq:theDiscriminant}
\Delta^\pm=(Tr^\pm) ^2 - 4Det^\pm = (\lambda + \epsilon(\alpha_{\epsilon0R} \pm \beta_{\epsilon0R}))(\lambda + \epsilon(\alpha_{\epsilon0R} 
\pm \beta_{\epsilon0R}) \mp 4\epsilon C_{det}) - 4\epsilon^2\beta^2_{\epsilon0I},
\end{equation}
we find that the eigenvalues of the 2x2 block of the Jacobian matrix \eqref{eq:jacobianMatrixSD} write as,
\begin{equation}\label{eq:transEigen}
 \begin{aligned}
 \bar{\mu}^{\pm}_2 &= -(\lambda + \epsilon(\alpha_{\epsilon0R} \pm 3\beta_{\epsilon0R})) - \sqrt{\xi},\\
  \bar{\mu}^{\pm}_3 &= -(\lambda + \epsilon(\alpha_{\epsilon0R} \pm 3\beta_{\epsilon0R})) + \sqrt{\xi},
 \end{aligned}
 \end{equation}
where 
$$
\xi=\big(\lambda + \epsilon(\alpha_{\epsilon0R} \pm \beta_{\epsilon0R})\big)\big(\lambda + \epsilon(\alpha_{\epsilon0R} \pm \beta_{\epsilon0R}) \mp 4\epsilon C_{det}\big) - 4\epsilon^2\beta^2_{\epsilon0I}.
$$

Next, we study the stability of the solutions $\bar{\mathcal{S}}^\pm_{osc}$ given in \eqref{eq:theSoscSols} when the parameters $\lambda,\epsilon$ lie in  
the domain
\begin{equation}\label{eq:areaA}
\mathcal{A}^{\pm}: = \Big\{(\lambda, \epsilon) \in \mathbb{R}^2 \quad | \quad \bar{\alpha}^\pm \geq 0, \quad \epsilon > 0 \Big\},
\end{equation}
where $\bar{\alpha}^\pm$ are defined in \eqref{eq:lesAlfes}. Notice that the domain $\mathcal{A}^{\pm}$ corresponds to the region on the right hand side of curves $C^{\pm}_{HB}$ and above the horizontal axis (see Fig. \ref{fig:nova2} left). Furthermore, as for the uncoupled case we link the solutions for the reduced system \eqref{eq:sdEqs} with the original system \eqref{eq:tereSystem}.

For $\bar{\alpha}^\pm = 0$, that is $(\lambda,\epsilon) \in C^{\pm}_{HB}$, the eigenvalues of the Jacobian matrix \eqref{eq:jacobianMatrixSD} at the fixed points $\bar{\mathcal{S}}^\pm_{osc}$ are given by
	\begin{equation}\label{eq:eigenvaluesHBCurve}
	\begin{aligned}
	\bar{\mu}^{\pm}_1 &= 0,\\
	\bar{\mu}^{\pm}_2 &= \mp 2\beta_{\epsilon0R} - i 2\epsilon\beta_{\epsilon0I},\\
	\bar{\mu}^{\pm}_3 &= \mp 2\beta_{\epsilon0R} + i 2\epsilon\beta_{\epsilon0I}.
	\end{aligned}
	\end{equation} 
	Therefore, when the parameters $(\lambda, \epsilon$) cross the curves $C_{HB}^\pm$ from left to right, if $\beta_{\epsilon0R} > 0$, $\bar{\mathcal{S}}^+_{osc}$ is a stable focus-node whereas $\bar{\mathcal{S}}^-_{osc}$ 
	is a saddle-focus with a $1$-dimensional stable manifold (corresponding to the $s$ direction which is always stable) and vice versa if $\beta_{\epsilon0R} < 0$. These results match exactly the results in Section \ref{sec:fixPointSec}: the 4D system has two periodic orbits that are born at different Hopf bifurcation curves $C^{+}_{HB}$ and $C^{-}_{HB}$ given in \eqref{eq:hbCurves}, and the stability of these periodic orbits depends on the sign of $\beta_{\epsilon 0R}$. 

For $\epsilon$ small and $\bar{\alpha}^\pm \geq 0$ the eigenvalues of the Jacobian matrix \eqref{eq:jacobianMatrixSD} at the fixed points $\bar{\mathcal{S}}^\pm_{osc}$ are given by
%
%
\begin{equation}\label{eq:eigenvaluesEpsZero}
\begin{aligned}
\bar{\mu}^{\pm}_1 &= -2\lambda + \mathcal{O}(\epsilon) ,\\
\bar{\mu}^{\pm}_2 &= -2\lambda + \mathcal{O}(\epsilon) ,\\
\bar{\mu}^{\pm}_3 &= \mp 2\epsilon (\beta_{\epsilon 0R} + C_{det}) + \mathcal{O}(\epsilon^2),
\end{aligned}
\end{equation}
which are $\mathcal{O}(\epsilon)$-close to the ones of the uncoupled case, $-2\lambda$ (double) and 0. 
In particular, depending on the sign of $(\beta_{\epsilon 0R} + C_{det})$, one fixed point is a stable node whereas the other is a saddle with a $1$-dimensional unstable manifold. 

\begin{figure}[H]
	\centering
	{\includegraphics[width=70mm]{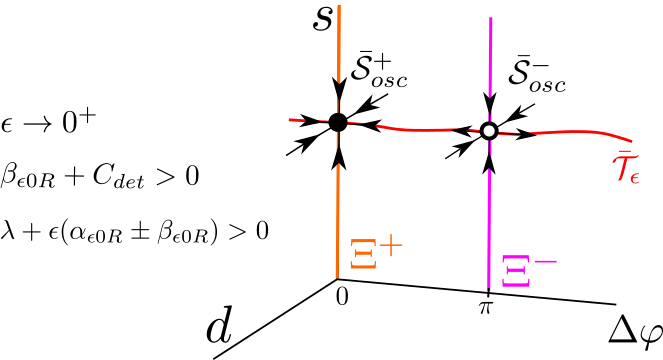}}
	\caption{Phase space of system \eqref{eq:sdEqs} for $\beta_{\epsilon0R} + C_{det} > 0$, $\lambda > 0$, and $0 < \epsilon < \epsilon_0(\lambda)$ (in particular   $\lambda + \epsilon(\alpha_{\epsilon0R} \pm \beta_{\epsilon0R}) > 0$). 
		There exist two fixed points $\bar{\mathcal{S}}^\pm_{osc}$, a stable node and a saddle point, respectively, which together with the unstable invariant manifold of the saddle point form the invariant curve $\bar {\mathcal{T}}_{\epsilon}$. 
		Due to the coupling term there are only two fixed points on $\bar{\mathcal{T}}_\epsilon$ whereas we had an infinite number in the unperturbed case. 
		Notice that the dynamics on the $s$ direction is always attracting.} \label{fig:unperturbedSD3}
\end{figure}

We remark that, for $\lambda>0$ fixed and $\epsilon$ small enough, we know that there exists an invariant curve $\bar{\mathcal{T}}_\epsilon$ corresponding to the 
invariant torus $\mathcal{T}_{\epsilon}$ obtained in Lemma~\ref{lem:persistence}.
Since this invariant curve is provided by Fenichel theory, it will contain the invariant points $\bar{\mathcal{S}}^{\pm}_{osc}$. Consequently, if $\beta_{\epsilon 0R}+C_{det}>0$, $\bar{\mathcal{T}}_\epsilon$ consists of the union of the saddle point 
$\bar{\mathcal{S}}^-_{osc}$, its unstable 1-dimensional manifold and the stable node
$\bar{\mathcal{S}}^+_{osc}$ (and vice versa if $\beta_{\epsilon0R}+C_{det} < 0$) (see Fig. \ref{fig:unperturbedSD3}). In conclusion, for $\lambda>0$ fixed and $\epsilon$ small enough, the invariant torus $\mathcal{T}_\epsilon$ of the system~\eqref{eq:tereSystem} contains the periodic orbits, $\mathcal{S}^+_{osc}$ and $\mathcal{S}^-_{osc}$ with $\Delta \varphi = 0$ and $\Delta \varphi = \pi$, respectively, 
whose stability depends on the sign of $\beta_{\epsilon 0R} + C_{det}$.

\begin{remark}\label{rem:rem2}
The existence of the invariant torus $\mathcal{T}_\epsilon$ is only guaranteed for $\lambda>0$ fixed and $\epsilon$ small enough by Lemma \ref{lem:persistence}. 
The evolution and eventual breakdown of this torus $\mathcal{T}_\epsilon$ (or, equivalently, the invariant curve $\bar{\mathcal{T}}_\epsilon$) when $\epsilon$ increases is beyond the scope of this paper.

However, in Section~\ref{sec:section4}, using system~\eqref{eq:sdEqs}, we study the evolution and bifurcations of the periodic orbits $\mathcal{S}^\pm_{osc}$ (corresponding to fixed points $\bar{\mathcal{S}}^\pm_{osc}$) for $(\lambda,\epsilon)$ small and no assumption on $\epsilon=o(\lambda)$. We leave as future work the exploration of the relationship between these bifurcations and the different mechanisms of destruction of the torus discussed in \cite{afraimovich1991invariant}.

\end{remark}

\section{Bifurcation diagrams of the oscillating $S_{osc}^{\pm}$ solutions}\label{sec:section4}

In the previous Sections we have shown that when $\epsilon$ is small and $\bar{\alpha}^\pm \geq 0$ there exist two critical points  $\bar{\mathcal{S}}^\pm_{osc}$ of system \eqref{eq:sdEqs} belonging to the 
curve $\bar{\mathcal{T}}_\epsilon$ which disappear at two independent curves: $C^\pm_{HB}$. Therefore, the points $\bar{\mathcal{S}}^\pm_{osc}$ undergo several bifurcations in the domain $\mathcal{A}^{\pm}$ defined in \eqref{eq:areaA}. Table \ref{tab:case1_table} shows the values of the trace $Tr^{\pm}$ in \eqref{eq:traceTereMatrix}, 
the determinant $Det^{\pm}$ in \eqref{eq:theDeterminant} and the discriminant $\Delta^{\pm}$ in \eqref{eq:theDiscriminant} of the Jacobian matrix of system \eqref{eq:sdEqs} at $\bar{\mathcal{S}}^\pm_{osc}$ near the curves $C^{\pm}_{HB}$ (given by the condition $\bar{\alpha}^\pm = 0$) and for $\bar{\alpha}^\pm \geq 0$ and $\epsilon$ small.
Notice that the sign of the constants $\beta_{\epsilon0R}$ and $C_{det} + \beta_{\epsilon0R}$ is relevant to determine the local dynamics around the fixed points. In particular,
\begin{itemize}
	\item  $\beta_{\epsilon0R}$ determines which of the two solutions $\bar{\mathcal{S}}^{\pm}_{osc}$ can have a null trace. 
	For $\beta_{\epsilon0R} > 0$, is $\bar{\mathcal{S}}^+_{osc}$, whereas for $\beta_{\epsilon0R} < 0$ is $\bar{\mathcal{S}}^-_{osc}$.

\item The sign of $C_{det} + \beta_{\epsilon0R}$ determines which of the two solutions $\bar{\mathcal{S}}^{\pm}_{osc}$ can have a null determinant. 
For $C_{det} + \beta_{\epsilon0R} > 0$, is $\bar{\mathcal{S}}^{-}_{osc}$, whereas for $C_{det} + \beta_{\epsilon0R} < 0$ is $\bar{\mathcal{S}}^{+}_{osc}$. 

\item Moreover, as we increase $\epsilon$ the discriminant always changes from negative to positive. That is, consistently with the eigenvalues obtained in \eqref{eq:eigenvaluesHBCurve} and \eqref{eq:eigenvaluesEpsZero}, 
the fixed points $\bar{\mathcal{S}}^{\pm}_{osc}$ change from a stable node and a saddle point to a stable focus and a saddle-focus.
\end{itemize}

Depending on the sign of $\beta_{\epsilon 0R}$ and $C_{det} + \beta_{\epsilon 0R}$ we consider three different cases: 
(1) $\beta_{\epsilon 0R}>0$, $C_{det} + \beta_{\epsilon 0R}>0$, (2) $\beta_{\epsilon 0R}<0$, $C_{det} + \beta_{\epsilon 0R}>0$, and (3) $\beta_{\epsilon 0R}=0$, $C_{det}>0$. The cases 
(i) $\beta_{\epsilon 0R}<0$, $C_{det} + \beta_{\epsilon 0R}<0$, (ii) $\beta_{\epsilon 0R}>0$, $C_{det} + \beta_{\epsilon 0R}<0$, and (iii) $\beta_{\epsilon 0R}=0$, $C_{det}<0$ are analogous to (1), (2) and (3), respectively, 
just replacing $\bar{\mathcal{S}}^{\pm}_{osc}$ by $\bar{\mathcal{S}}^{\mp}_{osc}$.
For each case, we study in detail the different bifurcations of the solutions $\bar{\mathcal{S}}^{\pm}_{osc}$ in the $(\lambda, \epsilon)$ parameter space, we link results obtained for the 3D system \eqref{eq:sdEqs} with the complete 4D system \eqref{eq:tereSystem}, and we discuss the regions of bistability.

\begin{table}[H]
	\renewcommand{\arraystretch}{1.75}
		\centering
		\scalebox{1}{
			\begin{tabular}{ c|c|c|c|c| }
				\cline{2-5}
				& \multicolumn{2}{ c| }{$\bar{\mathcal{S}}^+_{osc}$} & \multicolumn{2}{ c| }{$\bar{\mathcal{S}}^-_{osc}$} \\ \cline{2-5}
				& $\bar{\alpha}^+ \rightarrow 0^+$ &
				$\bar{\alpha}^+ \geq 0$, $\quad \epsilon \rightarrow 0^+$ & $\bar{\alpha}^- \rightarrow 0^+$ &
				$\bar{\alpha}^- \geq 0$, $\quad \epsilon \rightarrow 0^+$  \\ \hline
				\multicolumn{1}{|c|}{$Tr$} & $-4 \epsilon \beta_{\epsilon 0R}$ & $-2\lambda$ & $4 \epsilon \beta_{\epsilon 0R}$ & $-2\lambda$ \\ \hline
				\multicolumn{1}{|c|}{$Det$} & $4 \epsilon^2 (\beta^2_{\epsilon 0I} + \beta^2_{\epsilon 0R})$ & $4 \epsilon \lambda(C_{det}+\beta_{\epsilon 0R})$ & $4 \epsilon^2 (\beta^2_{\epsilon 0I} + \beta^2_{\epsilon 0R})$ & $-4 \epsilon \lambda(C_{det}+\beta_{\epsilon 0R})$ \\ \hline
				\multicolumn{1}{|c|}{$\Delta$} & $-4 \epsilon^2 \beta^2_{\epsilon 0I}$ & $ \lambda^2$ & $-4 \epsilon^2 \beta^2_{\epsilon 0I}$ & $ \lambda^2$ \\ \hline
			\end{tabular}}
			\renewcommand{\arraystretch}{1}
			\\
			\caption{Values for the trace ($Tr$), the determinant ($Det$) and the discriminant ($\Delta$) of the linearisation of system \eqref{eq:sdEqs} at the fixed points 
				$\bar{\mathcal{S}}^{\pm}_{osc}$  near the curves $C^{\pm}_{HB}$ ($\bar{\alpha}^\pm = 0$) and near to the uncoupled case ( $\bar{\alpha}^\pm \geq 0$ and $\epsilon$ small).} 
\label{tab:case1_table}
\end{table}

\subsection{Case $\beta_{\epsilon0R} > 0$ and $C_{det} + \beta_{\epsilon0R} > 0$  
\big( or $\beta_{\epsilon0R} < 0$ and $C_{det} + \beta_{\epsilon0R} < 0$ \big)}\label{sec:stableStab}

%

\subsubsection{Dynamics of $\bar{\mathcal{S}}^+_{osc}$}\label{sec:dynSoscCase1}

\begin{figure}[H]
	\centering
	{\includegraphics[width=80mm]{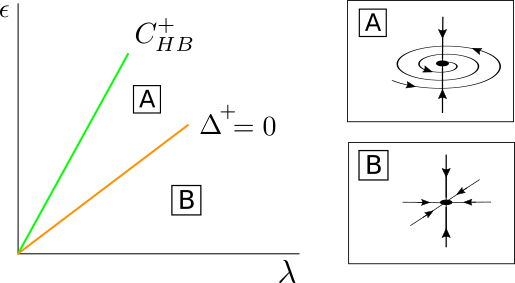}}
	\caption{Bifurcation diagram for $\bar{\mathcal{S}}^{+}_{osc}$ in the case $\beta_{\epsilon0R} > 0$ and $C_{det} + \beta_{\epsilon0R} > 0$. The fixed point $\bar{\mathcal{S}}^{+}_{osc}$ appears at 
	a supercritical pitchfork bifurcation of the origin occurring at the curve $C^+_{HB}$.} \label{fig:oscStable}
\end{figure}

For $\bar{\alpha}^+ \geq 0$, $\lambda$ fixed and $\epsilon$ small, the fixed point $\bar{\mathcal{S}}^{+}_{osc}$ for system \eqref{eq:sdEqs} is a stable node contained in the invariant curve $\bar{\mathcal{T}}_\epsilon$ 
(region B in Fig. \ref{fig:oscStable}), and as $\epsilon$ increases it becomes a stable focus at the curve $\Delta^+ = 0$ (region A in Fig.\ref{fig:oscStable}). It disappears at a pitchfork bifurcation of the origin in the $s$-direction at $C^+_{HB}$.

Going back to the original 4D system \eqref{eq:tereSystem} we have that for $\epsilon$ small there exists a stable periodic orbit $\mathcal{S}^+_{osc}$ 
(which belongs to the invariant torus $\mathcal{T}_\epsilon$), 
which disappears at a Hopf bifurcation of the origin  in $C^+_{HB}$. 

\subsubsection{Dynamics of $\bar{\mathcal{S}}^-_{osc}$}

The fixed point $\bar{\mathcal{S}}^-_{osc}$ changes from a saddle-focus with a $1$-dimensional stable manifold near $C^-_{HB}$ to a  saddle with a $2$-dimensional stable manifold for $\epsilon$ small and $\bar{\alpha}^- > 0$. 
Moreover, in this case the trace for $\bar{\mathcal{S}}^-_{osc}$ vanishes. Therefore, if 
\begin{equation}\label{eq:hopfInequality}
	\beta_{\epsilon0R} < - C_{det} + \sqrt{C^2_{det} + \beta^2_{\epsilon0I}},
\end{equation}
then $Tr^- = 0$ and $\Delta^- < 0$ and $\bar{\mathcal{S}}^-_{osc}$ undergoes a Hopf bifurcation.

So, we will distinguish two cases:

\subsubsection*{1) Case $\beta_{\epsilon0R} < - C_{det} + \sqrt{C^2_{det} + \beta^2_{\epsilon0I}},$}

For $\bar{\alpha}^+ \geq 0$, $\lambda$ fixed and $\epsilon$ small, the fixed point $\bar{\mathcal{S}}^-_{osc}$ is a saddle point with a $1$-dimensional 
unstable manifold (in the $\Delta \varphi$ direction) contained in the invariant curve $\bar{\mathcal{T}}_\epsilon$ (region D in Fig. \ref{fig:oscStable}). 
When crossing the curve $Det^-=0$ (region C), the point $\bar{\mathcal{S}}^-_{osc}$ becomes a stable node. 
As the coupling $\epsilon$ is increased, $\bar{\mathcal{S}}^{-}_{osc}$ crosses the curve $\Delta^-=0$ and $\bar{\mathcal{S}}^-_{osc}$ becomes a stable focus (region B).
When the parameters cross the curve $Tr^- = 0$, $\bar{\mathcal{S}}^-_{osc}$ undergoes a Hopf bifurcation $\bar{\mathcal{H}}$ in the $d, \Delta \varphi$ directions and 
$\bar{\mathcal{S}}^{-}_{osc}$ becomes a saddle focus with a $1$-dimensional unstable manifold (region A). 
At this bifurcation there appears or disappears a periodic orbit $\bar{\mathcal{T}}^-$ depending whether the Hopf bifurcation is supercritical or subcritical. 
Finally, the fixed point $\bar{\mathcal{S}}^{-}_{osc}$ disappears at a pitchfork bifurcation of the origin in the $s$-direction occurring at the curve $C^-_{HB}$.
\begin{figure}[H]
	\centering
	{\includegraphics[width=100mm]{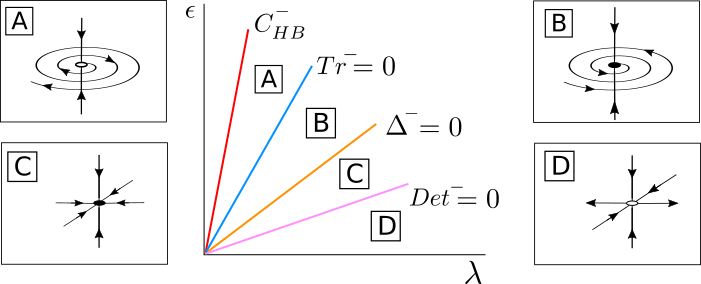}}
	\caption{Bifurcation diagram for $\bar{\mathcal{S}}^{-}_{osc}$ in the case $\beta_{\epsilon0R} > 0$, $C_{det} + \beta_{\epsilon0R} > 0$ and 
	$\beta_{\epsilon0R} < - C_{det} + \sqrt{C^2_{det} + \beta^2_{\epsilon0I}}$. The fixed point $\bar{\mathcal{S}}^{-}_{osc}$ appears at a supercritical pitchfork bifurcation of the origin occurring at the curve $C^-_{HB}$, 
	undergoes a Hopf bifurcation $\bar{\mathcal{H}}$ at the curve $Tr^-=0$ and becomes unstable at the curve $Det^-=0$.} \label{fig:oscUnstable_A}
\end{figure}
 
Going back to the original full 4D system \eqref{eq:tereSystem}, for $\epsilon$ small enough, there exists an unstable periodic orbit $\mathcal{S}^-_{osc}$, belonging to the torus 
${\mathcal{T}}_\epsilon$, which will become stable at the curve $Det^-=0$. The periodic orbit undergoes a Torus bifurcation and $\mathcal{S}^-_{osc}$ becomes unstable at the curve $Tr^-=0$ and a new torus $\mathcal{T}^-$ appears
or disappears depending whether the Torus bifurcation is subcritical or supercritical. 
Finally, $\mathcal{S}^-_{osc}$ will disappear at a Hopf bifurcation of the origin occurring at $C^-_{HB}$.

 \subsubsection*{2) Case $\beta_{\epsilon0R} > - C_{det} + \sqrt{C^2_{det} + \beta^2_{\epsilon0I}}$}

\begin{figure}[H]
	\centering
	{\includegraphics[width=80mm]{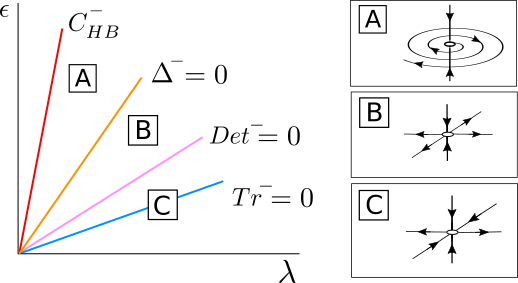}}
	\caption{Bifurcation diagram for $\bar{\mathcal{S}}^{-}_{osc}$ in the case $\beta_{\epsilon0R} > 0$, $C_{det} + \beta_{\epsilon0R} > 0$ and $\beta_{\epsilon0R} > - C_{det} + \sqrt{C^2_{det} + \beta^2_{\epsilon0I}}$. 
	The fixed point $\bar{\mathcal{S}}^{-}_{osc}$ appears at a supercritical pitchfork bifurcation of the origin in the $s$ direction occurring at the curve $C^-_{HB}$ and undergoes a bifurcation at the curve $Det^-=0$.} \label{fig:oscUnstable_B}
\end{figure}

For $\bar{\alpha}^+ \geq 0$, $\lambda$ fixed and $\epsilon$ small, the fixed point  
$\bar{\mathcal{S}}^-_{osc}$ is a saddle point with a $1$-dimensional unstable manifold 
(in the $\Delta \varphi$ direction) contained in  the invariant curve $\bar{\mathcal{T}}_\epsilon$ (region C in Fig. \ref{fig:oscUnstable_B}). As $\epsilon$ increases, $\bar{\mathcal{S}}^-_{osc}$ becomes a saddle with a $2$-dimensional 
unstable manifold at the curve $Det^-=0$ (region B). When further increasing the coupling $\epsilon$, $\bar{\mathcal{S}}^{-}_{osc}$ becomes a saddle-focus point at the curve $\Delta^-=0$ (region A), 
which disappears at a pitchfork bifurcation of the origin in the $s$-direction occurring at the curve $C^-_{HB}$.

Going back to the original full 4D system \eqref{eq:tereSystem}, for $\epsilon$ small enough, there exists an unstable periodic orbit $\mathcal{S}^-_{osc}$ belonging to the torus 
${\mathcal{T}}_\epsilon$. The periodic orbit undergoes a bifurcation at the curve $Det^-=0$ in which a stable manifold becomes unstable. Finally, $\mathcal{S}^-_{osc}$ will disappear at a Hopf bifurcation of the origin occurring at $C^-_{HB}$.

%


\subsubsection{Regions of bistability}

Since $\bar{\mathcal{S}}^+_{osc}$ is always stable, bistability between fixed points will appear in those regions where $\bar{\mathcal{S}}^-_{osc}$ is also stable. As in the case $\beta_{\epsilon0R} > - C_{det} + \sqrt{C^2_{det} + \beta^2_{\epsilon0I}}$, the
fixed point $\bar{\mathcal{S}}^-_{osc}$ is never stable, it is not possible to find bistability regions. By contrast, if $\beta_{\epsilon0R} < - C_{det} + \sqrt{C^2_{det} + \beta^2_{\epsilon0I}}$, there exist a region in the 
($\lambda, \epsilon$) parameter space defined as
\begin{equation}\label{eq:bistRegion1}
Tr^-(\lambda, \epsilon) < 0 \quad \text{and} \quad  Det^-(\lambda, \epsilon) > 0,
\end{equation}
in which $\bar{\mathcal{S}}^-_{osc}$ can be either a stable node or a stable focus (see Fig. \ref{fig:oscUnstable_A}). Thus, the system is bistable in the region \eqref{eq:bistRegion1}.

Moreover, in the case $\beta_{\epsilon0R} < - C_{det} + \sqrt{C^2_{det} + \beta^2_{\epsilon0I}}$, the point $\bar{\mathcal{S}}^-_{osc}$ undergoes a Hopf bifurcation $\bar{\mathcal{H}}$. If the Hopf bifurcation is supercritical, then $\bar{\mathcal{S}}^-_{osc}$ becomes unstable and a stable limit cycle $\bar{\mathcal{T}}^-$ appears, generating bistability between $\bar{\mathcal{S}}^+_{osc}$ and $\bar{\mathcal{T}}^-$. The detailed analysis of this situation is beyond 
the scope of this paper.

Finally, we remark that the same bistable scenarios can be found in the full system \eqref{eq:tereSystem} replacing the fixed points $\bar{\mathcal{S}}^\pm_{osc}$ by the limit cycles $\mathcal{S}^\pm_{osc}$ and 
the periodic orbit $\bar{\mathcal{T}}^-$ by the torus $\mathcal{T}^-$.

\subsection{Case $\beta_{\epsilon0R} < 0$ and $C_{det} + \beta_{\epsilon0R} > 0$  
\big(or $\beta_{\epsilon0R} > 0$ and $C_{det} + \beta_{\epsilon0R} < 0$ \big)}\label{sec:alterStb}


\subsubsection{Dynamics of $\bar{\mathcal{S}}^+_{osc}$}

In this case the trace for $\bar{\mathcal{S}}^+_{osc}$ vanishes ($Tr^+ = 0$). Therefore, as
\begin{equation}
\beta_{\epsilon0R} < - C_{det} < - C_{det} + \sqrt{C^2_{det} + \beta^2_{\epsilon0I}},
\end{equation}
then $Tr^+ = 0$ and $\Delta^+ < 0$ and $\bar{\mathcal{S}}^+_{osc}$ will always undergo a Hopf bifurcation $\bar{\mathcal{H}}$.

\begin{figure}[H]
	\centering
	{\includegraphics[width=80mm]{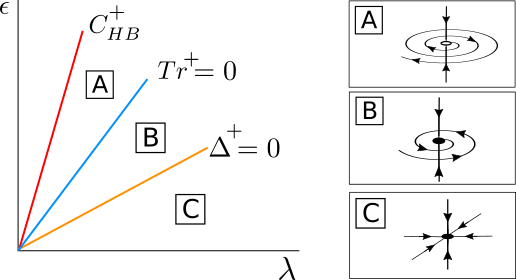}}
	\caption{Phase space for $\bar{\mathcal{S}}^{+}_{osc}$ in the case $\beta_{\epsilon0R} < 0$ and $C_{det} + \beta_{\epsilon0R} > 0$. The fixed point $\bar{\mathcal{S}}^{+}_{osc}$ appears at a supercritical pitchfork bifurcation of 
	the origin in the $s$ direction occurring at the curve $C^+_{HB}$, and undergoes a Hopf bifurcation $\bar{\mathcal{H}}$ at the curve $Tr^+=0$.} \label{fig:oscUnstableToStable}\end{figure}

For $\bar{\alpha}^+ \geq 0$, $\lambda$ fixed and $\epsilon$ small, the fixed point $\bar{\mathcal{S}}^+_{osc}$ is a stable node (region C in Fig. \ref{fig:oscUnstableToStable}) and becomes a stable focus when the parameters cross the 
curve $\Delta^+=0$ (region B). For larger values of $\epsilon$, the fixed point $\bar{\mathcal{S}}^+_{osc}$ undergoes a Hopf bifurcation $\bar{\mathcal{H}}$ at the curve $Tr^+=0$ and becomes a saddle-focus point (region A). 
At this bifurcation there appears or disappears a limit cycle $\bar{\mathcal{T}}^+$ depending whether this Hopf bifurcation is subcritical or supercritical. For larger values of $\epsilon$, the fixed point $\bar{\mathcal{S}}^+_{osc}$ 
disappears at a pitchfork bifurcation of the origin in the $s$-direction at the curve $C^+_{HB}$.

Going back to the original 4D system \eqref{eq:tereSystem}, for $\epsilon$ small enough there exists a stable periodic orbit $\mathcal{S}^+_{osc}$. This stable periodic orbit will lose its stability across a torus bifurcation occurring 
at the curve $Tr^+ = 0$. At this bifurcation there appears or disappears a torus $\mathcal{T}^+$ depending whether the torus bifurcation is subcritical or supercritical. Finally the unstable limit cycle $\mathcal{S}^+_{osc}$ collapses to the origin 
at a Hopf bifurcation occurring at the curve $C^+_{HB}$.

%
%

\subsubsection{Dynamics of $\bar{\mathcal{S}}^-_{osc}$}

\begin{figure}[H]
	\centering
	{\includegraphics[width=80mm]{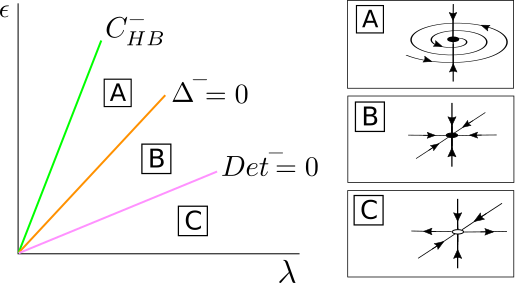}}
	\caption{Bifurcation diagram for $\bar{\mathcal{S}}^{-}_{osc}$ in the case $\beta_{\epsilon0R} < 0$ and $C_{det} + \beta_{\epsilon0R} > 0$. The fixed point $\bar{\mathcal{S}}^{-}_{osc}$ appers at a 
	supercritical pitchfork bifurcation of the origin in the $s$ direction occurring at the curve $C^-_{HB}$, and undergoes a bifurcation at $Det^-=0$.} \label{fig:oscStableToUnstable}
\end{figure}

For $\bar{\alpha}^- \geq 0$, $\lambda$ fixed and $\epsilon$ small, the fixed point $\bar{\mathcal{S}}^{-}_{osc}$ of system \eqref{eq:sdEqs} is a saddle point with a $1$-dimensional unstable manifold 
in the $\Delta \varphi$ direction contained in $\bar{\mathcal{T}}_\epsilon$ (region C in Fig. \ref{fig:oscStableToUnstable}), and as $\epsilon$ increases it becomes a stable node when $\epsilon$ crosses the curve $Det^- = 0$ (region B). 
For larger values of $\epsilon$, the fixed point $\bar{\mathcal{S}}^{-}_{osc}$ becomes a stable focus at the curve $\Delta^- = 0$ (region A) and 
disappears at a pitchfork bifurcation of the origin in the $s$ direction at the curve $C^{-}_{HB}$.

Going back to the original 4D system \eqref{eq:tereSystem}, for $\epsilon$ small there exists an unstable periodic orbit $\mathcal{S}^-_{osc}$. 
This unstable periodic orbit becomes stable at the curve $Det^- = 0$. Finally, the stable limit cycle $\mathcal{S}^-_{osc}$ collapses to the origin at a Hopf bifurcation occurring at the curve $C^-_{HB}$.

\subsubsection{Regions of bistability}

There exist a region in the ($\lambda, \epsilon$)-parameter space given by
\begin{equation}
Tr^+(\lambda, \epsilon) < 0 \quad \text{and} \quad Det^-(\lambda, \epsilon)> 0,
\end{equation}
in which both fixed points $\bar{\mathcal{S}}^\pm_{osc}$ are stable. If the Hopf bifurcation is supercritical, then $\bar{\mathcal{S}}^+_{osc}$ becomes unstable and a stable limit cycle $\bar{\mathcal{T}}^+$ appears, generating bistability between $\bar{\mathcal{S}}^-_{osc}$ and $\bar{\mathcal{T}}^+$. The detailed analysis of this situation is beyond of the scope of this paper.

Finally, we remark that the same bistable scenarios can be found in the full system \eqref{eq:tereSystem} replacing the fixed points $\bar{\mathcal{S}}^\pm_{osc}$ by the limit cycles $\mathcal{S}^\pm_{osc}$ and the periodic orbit $\bar{\mathcal{T}}^+$ by the torus $\mathcal{T}^+$.

\subsection{The degenerated case $\beta_{\epsilon0R} = 0$ and $C_{det} > 0$ \big(or $\beta_{\epsilon0R} = 0$ and $C_{det} < 0$ \big)}\label{sec:extDeg}

In this case, the curves $C^{\pm}_{HB}$ coincide. Moreover, the trace in \eqref{eq:traceTereMatrix} is identically zero for ($\lambda, \epsilon) \in C^{\pm}_{HB}$. 
To obtain the sign of $Tr^{\pm}$, we compute $Tr^{\pm}$ when $\lambda + \epsilon \alpha_{\epsilon0R} \rightarrow 0^+$. We have
\begin{equation}\label{eq:lastEquation}
Tr(\lambda, \epsilon) =  (\lambda + \epsilon \alpha_{\epsilon0R}) \big( -2 +  \mathcal{O}_2(\epsilon) \big),
\end{equation}
so, near the $C_{HB}^\pm$ curves, both fixed points $\bar{\mathcal{S}}^{\pm}_{osc}$ are stable. 

\subsubsection{Dynamics of $\bar{\mathcal{S}}^+_{osc}$}

\begin{figure}[H]
	\centering
	{\includegraphics[width=80mm]{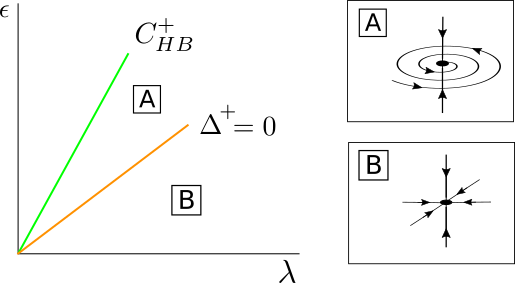}}
	\caption{Bifurcation diagram for $\bar{\mathcal{S}}^{+}_{osc}$ in the case $\beta_{\epsilon0R} = 0$ and $C_{det} > 0$. The fixed point $\bar{\mathcal{S}}^{+}_{osc}$ apperas at a supercritical pitchfork bifurcation 
	of the origin in the $s$ direction occurring at the curve $C^+_{HB}$.} \label{fig:oscExtremelyStable}
	\end{figure}

For $\bar{\alpha}^+ \geq 0$, $\lambda$ fixed and $\epsilon$ small,  the fixed point $\bar{\mathcal{S}}^{+}_{osc}$ is a stable node (region B in Fig. \ref{fig:oscExtremelyStable}), and as $\epsilon$ increases 
it becomes a stable focus when the parameters cross the curve $\Delta^+ = 0$ (region A). For larger values of $\epsilon$, the fixed point $\bar{\mathcal{S}}^{+}_{osc}$ disappears at a pitchfork bifurcation of the origin in the $s$ direction at the curve $C^+_{HB}$.

Going back to the original 4D system \eqref{eq:tereSystem}, for $\epsilon$ small there exists a stable periodic orbit  $\mathcal{S}^+_{osc}$, 
which collapses to the origin at a Hopf bifurcation occurring at the curve $C^+_{HB}$.

\subsubsection{Dynamics of $\bar{\mathcal{S}}^-_{osc}$}

\begin{figure}[H]
	\centering
	{\includegraphics[width=80mm]{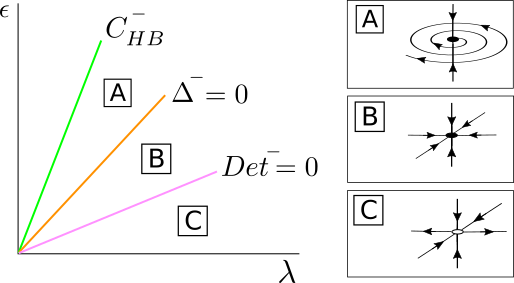}}
	\caption{Phase space for the $\bar{\mathcal{S}}^{-}_{osc}$ fixed point in the  case $\beta_{\epsilon0R} = 0$ and $C_{det} > 0$. The fixed point $\bar{\mathcal{S}}^{-}_{osc}$ apperas at a supercritical pitchfork bifurcation 
	of the origin in the $s$ direction occurring at the curve $C^-_{HB}$, and undergoes a bifurcation at the curve $Det^-=0$.} \label{fig:oscExtremelyUnstable}
\end{figure}

For $\bar{\alpha}^+ \geq 0$, $\lambda$ fixed and $\epsilon$ small, the fixed point $\bar{\mathcal{S}}^{-}_{osc}$ is a saddle point with a 1-dimensional unstable manifold (region C in Fig. \ref{fig:oscExtremelyUnstable}), and as $\epsilon$ increases 
it becomes a stable node when the parameters cross the curve $Det^- = 0$ (region B). For larger $\epsilon$ values the fixed point  $\bar{\mathcal{S}}^{-}_{osc}$ becomes a stable focus at $\Delta^-=0$ (region A) which collapses at a pitchfork bifurcation of the origin in the $s$ direction at the curve $C^-_{HB}$. 

Going back to the original 4D system \eqref{eq:tereSystem}, for $\epsilon$ small there exists an unstable periodic orbit $\mathcal{S}^-_{osc}$  which changes stability at the curve $Det^-=0$. 
Finally, the stable periodic orbit $\mathcal{S}^-_{osc}$ collapses to the origin at a Hopf bifurcation at the curve $C^-_{HB}$.  

%
%

\subsubsection{Regions of bistability} 

In the region in the ($\lambda, \epsilon$)-parameter space given by
\begin{equation}
	 Det^-(\lambda, \epsilon) > 0
\end{equation}
both fixed points $\bar{\mathcal{S}}^+_{osc}$ and $\bar{\mathcal{S}}^-_{osc}$ are stable.

We remark that the same bistability scenarios can be found in the full system \eqref{eq:tereSystem} replacing the fixed points $\bar{\mathcal{S}}^\pm_{osc}$ 
by the limit cycles $\mathcal{S}^\pm_{osc}$.

\section{Wilson-Cowan models for perceptual bistability}\label{sec:percept}

Wilson-Cowan oscillators are biophysically
motivated neural oscillators providing a population-averaged firing rate
description of neural activity, which have been widely used to study
cortical dynamics and cortical oscillations
\cite{Wilson1972,rinzel1998analysis}. The Wilson-Cowan oscillator
(an excitatory ($E$), inhibitory ($I$) pair) considered here has dynamics described
by
\begin{equation}\label{eq:oneWC}
\begin{split}
\tau \dot{E} &= -E + S(aE - bI),\\
\tau \dot{I} &= -I + S(cE - dI),
\end{split}
\end{equation}
where $\tau$ is a time constant and the constants $a,b,c$ and $d$ are the intrinsic $E$ to $E$, $I$ to $E$, $E$ to $I$ and $I$ to $I$
coupling weights, respectively. The function $S$ is the sigmoidal response function
\begin{equation}\label{eq:sigmoid}
S(x) = \frac{1}{1 + e^{-\lambda x + \theta}} - \frac{1}{1 + e^{\theta}},
\end{equation}
which has threshold $\theta$ and slope $\lambda$ with the convenient property $S(0) = 0$. 
The function $S$ has the property $S'(0) = \lambda S_1$, where $S_1 = \frac{e^{\theta}}{(1 + e^{\theta})^2}$, and $\lambda$ is treated as a 
bifurcation parameter playing the equivalent role to $\lambda$ in previous sections.

The system generically has a steady state
$(E,I)=(0,0)$, which undergoes a Hopf bifurcation at $\lambda_c=\frac{2}{(a-d)S_1}$. When
coupled with a second, identical oscillator the $4$-dimensional pair of Wilson-Cowan oscillators (E-I pairs) coupled with strength $\epsilon$ are given by
\begin{equation}
\begin{split}\label{eq:twoWCs}
\tau \dot{E}_1 &= -E_1 + S(aE_1 - bI_1),\\
\tau \dot{I}_1 &= -I_1 + S(cE_1 - dI_1 + \epsilon (E_2 - b_{sp}I_2)),\\
\tau \dot{E}_2 &= -E_2 + S(aE_2 - bI_2),\\
\tau \dot{I}_2 &= -I_2 + S(cE_2 - dI_2 + \epsilon (E_1 - b_{sp}I_1)),
\end{split}
\end{equation}
whose dynamics will be explored in this Section. 

For this study, we will consider the following set of parameters:
\begin{equation}
\mathcal{P} = \{  a=7, b=5.25, \\c=5, d=0.7, \theta=2, \tau=1 \},
\end{equation}
whereas $\lambda$ and $\epsilon$ will be the bifurcation parameters. 
By considering $b_{sp} = -0.03, 0.03, 0.0$ we will study different types of dynamics. 
For each case we will write system \eqref{eq:twoWCs} in the normal form \eqref{eq:tereSystem} by numerically computing its corresponding 
coefficients (see Appendix \ref{sec:normalFormCoefs}). 
Next, by using numerical continuation we will compute bifurcation diagrams for system \eqref{eq:twoWCs}, 
so we can check the theoretical predictions in Section \ref{sec:section4} and complete the bifurcation diagrams for large values of 
$\lambda$ and $\epsilon$, where the normal form approximation breaks down.

\begin{table}[H]
	\renewcommand{\arraystretch}{1.75}
	\begin{minipage}{.5\textwidth}
		\centering
		\scalebox{1.25}{
			\begin{tabular}{ |c|c|c|c| }
				\cline{2-4}
				\multicolumn{1}{ c| }{} & \multicolumn{3}{ c| }{$b_{sp}$} \\ \cline{2-4}
				\multicolumn{1}{ c| }{} & -0.03 & 0.03 & 0 \\ \hline
				$\alpha_{01R}$ & -21.94 & -21.94 & -21.94 \\ \hline
				$\alpha_{01I}$ & -20.94 & -20.94 & -20.94 \\ \hline
				$\alpha_{\epsilon 0R}$ & 0 & 0 & 0 \\ \hline
				$\alpha_{\epsilon 0I}$ & 0 & 0 & 0 \\ \hline
				$\alpha_{\epsilon 1R}$ & 0 & 0 & 0 \\ \hline
				$\alpha_{\epsilon 1I}$ & 0 & 0 & 0 \\ \hline
				$\alpha_{\epsilon 2R}$ & 8.4 & 9.02 & 8.72 \\ \hline
				$\alpha_{\epsilon 2I}$ & 6.34 & 6.8 & 6.57 \\ \hline
				$\alpha_{\epsilon 3R}$ & -24.02 & -22.3 & -23.2 \\ \hline
				$\alpha_{\epsilon 3I}$ & -46.36 & -44.92 & -45.46 \\ \hline
			\end{tabular}}
		\end{minipage}
		\begin{minipage}{.5\textwidth}
			\centering
			\scalebox{1.25}{
				\begin{tabular}{ |c|c|c|c| }
				\cline{2-4}
				\multicolumn{1}{ c| }{} & \multicolumn{3}{ c| }{$b_{sp}$} \\ \cline{2-4}
				\multicolumn{1}{ c| }{} & -0.03 & 0.03 & 0 \\ \hline
				$\omega$ & 1.073 & 1.073 & 1.073 \\ \hline
				$\beta_{\epsilon 0R}$ & 0.0047 & -0.0047 & 0 \\ \hline
				$\beta_{\epsilon 0I}$ & 0.252 & 0.241 & 0.246 \\ \hline
				$\beta_{\epsilon 1R}$ & -12.91 & -13.18 & -13.05 \\ \hline
				$\beta_{\epsilon 1I}$ & 19.36 & 16.76 & 18.06 \\ \hline
				$\beta_{\epsilon 2R}$ & 7.16 & 6.46 & 6.52 \\ \hline
				$\beta_{\epsilon 2I}$ & -5.56 & -5.47 & -5.52 \\ \hline
				$\beta_{\epsilon 3R}$ & 14.29 & 13.33 & 13.81 \\ \hline
				$\beta_{\epsilon 3I}$ & 10.02 & 10.3 & 10.16 \\ \hline
				\multicolumn{1}{ c }{}
				\end{tabular}}
				
			\end{minipage}
			\renewcommand{\arraystretch}{1}
			\\
			\caption{Coefficients of the normal form \eqref{eq:tereSystem} for the three considered cases, namely $b_{sp} = -0.03, 0.03$ and $0$. 
			These coefficients have been computed using the procedure described in Appendix \ref{sec:normalFormCoefs}.} 
			\label{tab:coefsTable}
		\end{table}


\subsection{Case $b_{sp} < 0$}

\begin{figure}[t]
	\centering
	{\includegraphics[width=120mm]{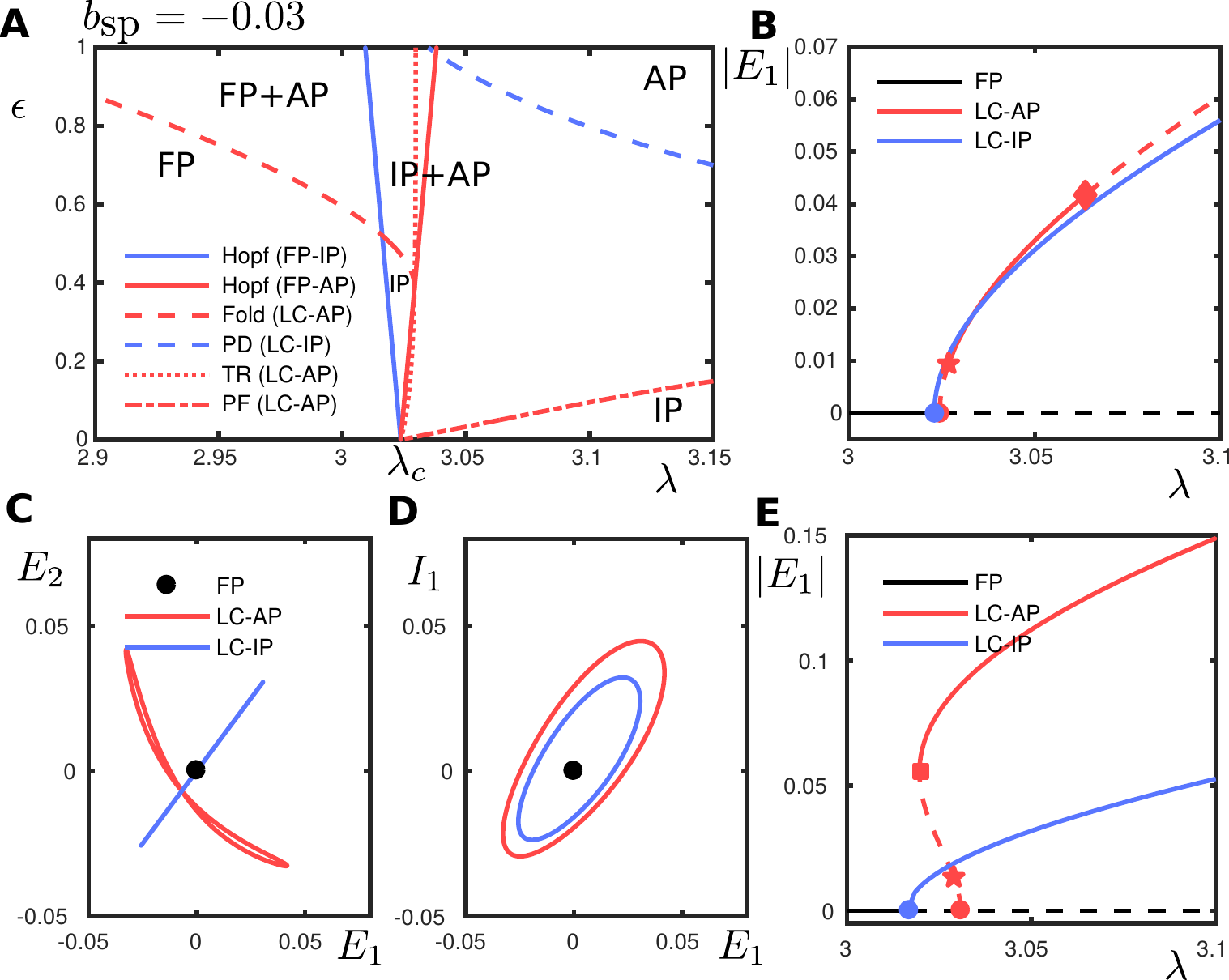}}
	\caption{Bifurcation diagram with parameters $\mathcal{P}$ and 
		$b_{sp}=-0.03$ in \eqref{eq:twoWCs} (corresponding to the case
		$\beta_{\epsilon0R} > 0$, $C_{det} + \beta_{\epsilon0R} > 0$ and satisfying
		$\beta_{\epsilon0R} < - C_{det} + \sqrt{C^2_{det} +
			\beta^2_{\epsilon0I}}$
		as described in Section \ref{sec:stableStab}). \textbf{A}:
		Two-parameter bifurcation diagram in the
		$(\lambda,\epsilon)$-plane.  The legend indicates
		bifurcations of a fixed point (FP) or a limit cycle (LC) giving rise
		to or involving the $\Delta \varphi = 0$ in-phase (IP) or $\Delta \varphi = \pi$ anti-phase (AP) solution branches; PD: period doubling; PF: pitchfork; TR: torus bifurcation. Text labels indicate
		the solutions that are stable in a given region, e.g. `IP+AP' is a
		region with coexisting, stable IP and AP solutions. \textbf{B}:
		One-parameter bifurcation diagram at $\eps=0.05$ showing the FP branch, IP branch
		and AP branch; dashed segments are unstable.  The IP and AP branches
		bifurcate from the FP branch in subsequent Hopf bifurcations
		(bullet) for $\lambda$ increasing.  The IP branch emerges stable and
		remains stable. For increasing $\lambda$ the AP branch is initially
		unstable, gains stability at a torus bifurcation (star) and loses
		stability at a pitchfork bifurcation
		(diamond). \textbf{C}: Coexisting solutions at $\lambda\approx3.05$ and $\epsilon=0.05$ in the $(E_1,E_2)$-plane. Motion on the diagonal (blue) corresponds to in-phase oscillations. \textbf{D}: As \textbf{C} in the $(E_1,I_1)$-plane for one E-I oscillator. \textbf{E}: As \textbf{C} at $\epsilon=0.5$, where a torus bifurcation (star) is on an unstable branch that gains stability at a Fold of limit cycle (square). }\label{fig:bifThirdCase}
\end{figure}

We consider the case $b_{sp} = -0.03$. 
The coefficients of the normal form, which were computed using the techniques described in Appendix \ref{sec:normalFormCoefs}, 
are given in Table \ref{tab:coefsTable} and satisfy the conditions $\beta_{\epsilon0R} > 0$, 
$C_{det} + \beta_{\epsilon0R} > 0$ and $\beta_{\epsilon0R} < - C_{det} + \sqrt{C^2_{det} + \beta^2_{\epsilon0I}}$. 
Therefore, this case corresponds to the one considered in Section \ref{sec:stableStab}. Fig. \ref{fig:bifThirdCase} shows the bifurcation diagram of 
system \eqref{eq:twoWCs} for $b_{sp} = -0.03$ obtained numerically. 
The results match the theoretical predictions obtained in Section \ref{sec:stableStab}. 
More precisely, for a fixed $\epsilon$ value and varying the bifurcation parameter $\lambda$ we have:
\begin{itemize}
	\item A stable in-phase (IP) solution corresponding to $\mathcal{S}^+_{osc}$ will emerge from the Hopf bifurcation at $C^+_{HB}$. 
	Moreover when varying the bifurcation parameter, the IP solution will maintain its stability (see Fig. \ref{fig:oscStable}).
	\item An unstable anti-phase (AP) solution corresponding to $\mathcal{S}^-_{osc}$ will emerge from the Hopf bifurcation at $C^-_{HB}$. 
	For a fixed $\epsilon$ and varying the bifurcation parameter AP solution gains stability across a Torus bifurcation, but when further increasing 
	the bifurcation parameter it will loose it again across a pitchfork bifurcation (corresponding respectively to the lines $Tr^-=0$ and $Det^-=0$ 
	in Fig. \ref{fig:oscUnstable_A}).
\end{itemize}

\subsection{Case $b_{sp} > 0$}

\begin{figure}[t]
	\centering
{\includegraphics[width=120mm]{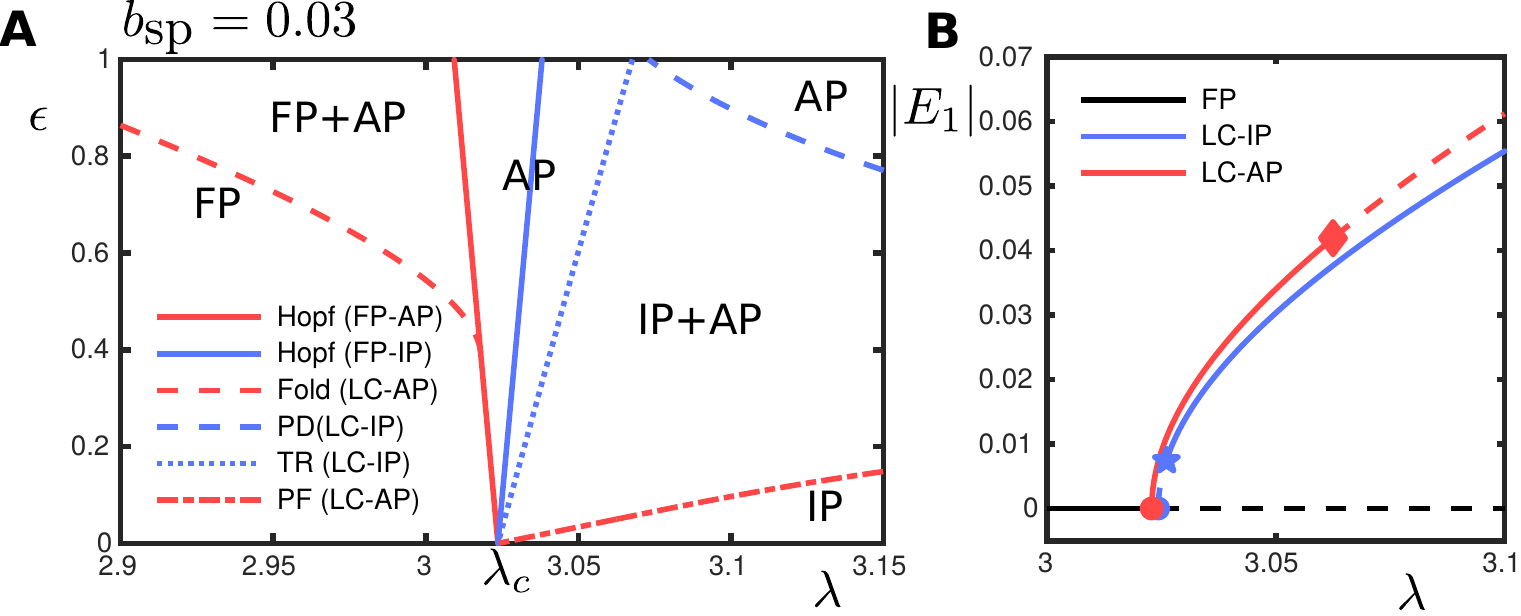}}
\caption{Bifurcation diagram with parameters $\mathcal{P}$ and
  $b_{sp}=0.03$ in \eqref{eq:twoWCs} (corresponding to the case
  $\beta_{\epsilon0R} < 0$ and $C_{det} + \beta_{\epsilon0R} > 0$, as
  described in Section \ref{sec:alterStb}). \textbf{A}: Two-parameter
  bifurcation diagram in the $(\lambda,\epsilon)$-plane. Legends and
  labelling as in Fig.  \ref{fig:bifThirdCase}; TR: torus
  bifurcation. \textbf{B}: One-parameter bifurcation diagram at $\eps=0.05$ showing
  the FP branch, IP branch and AP branch; dashed segments are
  unstable. The AP and IP branches bifurcate from the FP branch in
  subsequent Hopf bifurcations (bullet) for $\lambda$ increasing. The AP
  branch loses stability in a pitchfork bifurcation (diamond). The IP
  branch is initially unstable and gains stability at a torus
  bifurcation (star).}\label{fig:bifSecCase}
\end{figure} 

We consider the case $b_{sp} = 0.03$. 
The coefficients of the normal form, which were computed using the techniques described in Appendix \ref{sec:normalFormCoefs}, 
are given in Table \ref{tab:coefsTable} and satisfy the conditions $\beta_{\epsilon0R} < 0$ and $C_{det} + \beta_{\epsilon0R} > 0$. 
Therefore, this case corresponds to the one considered in Section \ref{sec:alterStb}. 
Fig. \ref{fig:bifSecCase} shows the bifurcation diagram of system \eqref{eq:twoWCs} for $b_{sp} = 0.03$ obtained numerically. 
The results match the theoretical predictions in Section \ref{sec:alterStb}. 
More precisely, for a fixed $\epsilon$ value and varying the bifurcation parameter $\lambda$ we have:
\begin{itemize}
\item A stable anti-phase (AP) solution corresponding to $\mathcal{S}^-_{osc}$ will emerge from a Hopf bifurcation at $C^{-}_{HB}$ whereas an unstable in-phase (IP) solution corresponding to $\mathcal{S}^+_{osc}$ will emerge from the Hopf bifurcation at $C^{+}_{HB}$
\item The stability of both solutions is reversed as the bifurcation parameter grows. Moreover, the bifurcations giving rise to these stability changes are of the same type as we predicted: IP solution becomes stable across a torus bifurcation (corresponding to the Hopf bifurcation $\bar{\mathcal{H}}$ at the $Tr^+=0$ line in Fig. \ref{fig:oscUnstableToStable}) whereas the AP solution looses stability across a pitchfork bifurcation of limit cycles (corresponding to the $Det^-=0$ line in Fig. \ref{fig:oscStableToUnstable}).
\end{itemize}


\subsection{Case $b_{sp} = 0$}

We consider the case $b_{sp} = 0.0$. The coefficients of the normal form, which were computed using the techniques described in Appendix \ref{sec:normalFormCoefs}, are given in Table \ref{tab:coefsTable} and satisfy the conditions $\beta_{\epsilon0R} = 0$ and $C_{det} > 0$. Therefore, this case corresponds to the ``degenerated case'' discussed in Section \ref{sec:extDeg}. Fig. \ref{fig:bifFirstCase} shows the bifurcation diagram of system \eqref{eq:twoWCs} for $b_{sp} = 0$ obtained numerically. Notice that it matches the theoretical predictions, namely: 

\begin{figure}[t]
	\centering
	{\includegraphics[width=120mm]{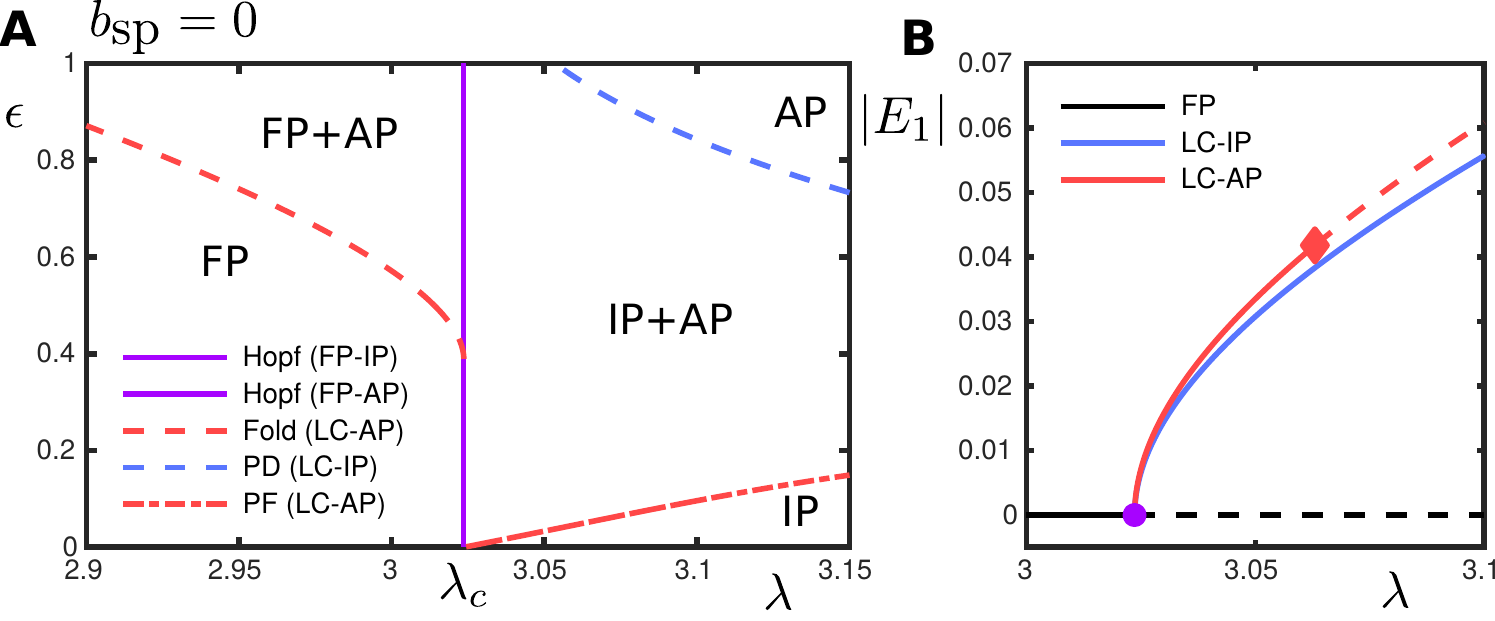}}
	\caption{Bifurcation diagram with parameters $\mathcal{P}$ and
		$b_{sp}=0$ in \eqref{eq:twoWCs} (corresponding to the ``Degenerated case'' in section \ref{sec:extDeg}). \textbf{A}:
		Two-parameter bifurcation diagram where curves are the locus of
		bifurcations in the $(\mu,\epsilon)$-plane. The legend indicates
		bifurcations of a fixed point (FP) or a limit cycle (LC) giving rise
		to or involving the $\Delta \varphi = 0$ in-phase (IP) or $\Delta \varphi = \pi$ anti-phase (AP) solution branches; PD: period doubling; PF: pitchfork. Text labels indicate
		the solutions that are stable in a given region, e.g. `IP+AP' is a
		region with coexisting, stable IP and AP solutions. \textbf{B}:
		One-parameter bifurcation diagram for fixed $\epsilon=0.05$ showing
		the fixed point branch, IP branch and AP branch; dashed segments are
		unstable. The IP and AP branches bifurcation from the FP branch at a
		degenerate Hopf bifurcation (bullet). The AP branch loses stability
		in a pitchfork bifurcation (diamond). } \label{fig:bifFirstCase}
\end{figure}

\begin{itemize}
	\item Both Hopf bifurcation curves $C^{\pm}_{HB}$ coincide and give rise to a bistable situation. 
	On one side of the double Hopf curve there exists bistability between the in-phase (IP) solution $\Delta \varphi = 0$ corresponding to 
	$\mathcal{S}^+_{osc}$ and the anti-phase (AP) $\Delta \varphi = \pi$ solution corresponding to $\mathcal{S}^-_{osc}$.
	\item For $\epsilon$ fixed and increasing the bifurcation parameter $\lambda$, the $\mathcal{S}^-_{osc}$ (AP) solution loses stability across a 
	Pitchfork bifurcation of limit cycles that we found for the 3D system as the line having $Det^-=0$ (see Fig. \ref{fig:oscExtremelyUnstable}).
\end{itemize}

\subsection{Dynamics beyond the weak coupling limit}

Our numerical bifurcation analysis has revealed the possibility for richer dynamics, whilst noting a wide range of parameters for which the IP and AP 
solutions are stable and coexist. 
Furthermore, a Bautin bifurcation on the AP Hopf branch for $\epsilon_{BT}\approx0.4$  as seen in 
Figures~\ref{fig:bifThirdCase},~\ref{fig:bifSecCase} and~\ref{fig:bifFirstCase} gives rise to a region of parameter space for 
$\lambda\lesssim\lambda_c$ where a stable AP solution coexists with a stable fixed point. The bifurcation point $\epsilon_{BT}$ separates branches of sub- and supercritical Hopf bifurcations in the parameter space. As we can see, for nearby $\lambda, \epsilon$ parameter values, the system has two limit cycles which collide and disappear via a Fold bifurcation of periodic orbits. Although the analysis done in Sections \ref{sec:section3} and \ref{sec:section4} is restricted to the weak coupling case, we briefly discuss how the reduced system \eqref{eq:sdEqs} can provide some insight about this bifurcation.


In the weak coupling regime, the denominator in the formula \eqref{eq:equivariantSolutions} for the $s^\pm_{osc}$ solutions, is given by $\alpha_{01R} + \epsilon K_{stb}^\pm$  and is assumed to be negative. Therefore, $s^\pm_{osc}$ solutions appear for $\alpha^\pm = \lambda + \epsilon(\alpha_{\epsilon 0R} \pm \beta_{\epsilon 0R}) > 0$ at a supercritical pitchfork bifurcation of the origin (see Fig. \ref{fig:perturbedSDpitch}). Nevertheless, writing the equation for $s$ in \eqref{eq:equivariantSub} in the following way
\[\dot{s}=A(\lambda,\epsilon)s + B(\lambda,\epsilon)s^3,\]
we clearly see that at the curve $A(\lambda,\epsilon)=0$, the origin undergoes a pitchfork bifurcation that it supercritical or subcritical depending on the sign of 
$B(\lambda,\epsilon)$. Consequently, the point ($\lambda, \epsilon$) satisfying $A(\lambda,\epsilon) = 0$ and $B(\lambda,\epsilon) = 0$ corresponds to a Bautin bifurcation. Thus, using the expression for $A$ and $B$ (which are known up to first order in $\epsilon$ and $\lambda$), we can estimate that a Bautin bifurcation occurs for 
\begin{equation}
\epsilon_{BT} \approx -\frac{\alpha_{01R}}{K^-_{stb}},
\end{equation}
assuming that $K^-_{stb}>0$ and for $\lambda_{BT}$ such that $(\lambda_{BT}, \epsilon_{BT}) \in C_{HB}^-$. Although an accurate derivation is beyond the scope of this work, this transition from subcritical to supercritical involves the appearance of a curve of saddle-node bifurcations of fixed points for the system \eqref{eq:sdEqs} for nearby values of the parameters. More precisely, if we consider the exact expression of the determinant of the 2x2 block of Jacobian Matrix \eqref{eq:jacobianMatrixSD} given by:
\begin{equation}
Det(s_{osc}^-) = c_d^d c_\varphi^\varphi - c_d^\varphi c_\varphi^d 
\end{equation}
where the constants are given by Eqs. \eqref{eq:jacobianTermsCs} in the Appendix \ref{sec:appendix} with $s = s_{osc}^-$ in \eqref{eq:equivariantSolutions}, one can see that it is singular at $B(\lambda,\epsilon) = 0$. Therefore, we consider the curve 
\begin{equation*}
B(\lambda,\epsilon) Det(s_{osc}^-) = 0,
\end{equation*}
and one can see that the Bautin point ($\lambda_{BT}, \epsilon_{BT}$) belongs to it. Moreover, for $\epsilon > \epsilon_{BT}$ as $B(\lambda, \epsilon) > 0$ this curve corresponds to the saddle-node bifurcations of the solutions $s_{osc}^{-}$ outside the $C_{HB}^-$ curve.

Using the numerical values given in Table~\ref{tab:coefsTable}, $K_{stb}^->0$. Thus, we can estimate from the normal form that the Bautin bifurcation occurs for $\epsilon_{BT} \approx 0.42, 0.43, 0.42$ for $b_{sp} = -0.03, 0.03, 0$, respectively, which matches the results obtained numerically (see Figures~\ref{fig:bifThirdCase}, \ref{fig:bifSecCase} and~\ref{fig:bifFirstCase}).
Recall that in the original 4D system \eqref{eq:tereSystem} the pitchfork and saddle-node bifurcations correspond to Hopf and fold of limit cycles bifurcations, respectively.

Besides this previous behaviour, we also remark that the IP solution undergoes a period-doubling bifurcation for large $\epsilon$ and $\lambda$ leading to richer dynamical behaviour 
away from the analytically-investigated uncoupling limit.

\subsection{Periodically Forced Coupled Wilson-Cowan Equations}

With the aim of finding coexisting IP and AP solutions (corresponding
to ``percept 1'' and ``percept 2'' as described in section
\ref{sec:intro}) we now introduce periodic forcing terms to the
coupled WC system given by \eqref{eq:twoWCs}. We consider anti-phase
inputs with forcing frequency $f=2.5\,Hz$ and amplitude $A$ which will
be varied as a bifurcation parameter: 
\begin{equation}\label{eq:twoWCPF}
\begin{split}
\tau \dot{E}_1 &= -E_1 + S(aE_1 - bI_1+A\sin^{2n}(2\pi f t)+(1-h)A\cos^{2n}(2\pi f t)),\\
\tau \dot{I}_1 &= -I_1 + S(cE_1 - dI_1 + \epsilon (E_2 - b_{sp}I_2)),\\
\tau \dot{E}_2 &= -E_2 + S(aE_2 - bI_2+A\cos^{2n}(2\pi f t)+(1-h)A\sin^{2n}(2\pi f t))),\\
\tau \dot{I}_2 &= -I_2 + S(cE_2 - dI_2 + \epsilon (E_1 - b_{sp}I_1)),
\end{split}
\end{equation}
where the parameters $\mathcal{P}$ (with the exception of $\tau$) and
nonlinearity \eqref{eq:sigmoid} are as above. The input asymmetry
parameter $h$ controls the balance of inputs across the two
oscillators; when $h=1$ the oscillators receive exclusive inputs (the
case typically considered in competition models
\cite{Laing2002,Wilson2003,Shpiro2009,li2017attention}) and when $h=0$
the oscillators receive identical inputs (the case considered
here). The forcing terms are raised to an even power $2n$ with $n=5$
to be positive and \emph{sharpened} such that the anti-phase inputs do
not overlap in time. Noting that the isolated Wilson-Cowan
oscillator \label{eq:twoWC} undergoes a supercritical Hopf bifurcation
at $\lambda=\frac{2}{(a-d)S_1}=3.025$, we set $\lambda=2.6$, before this
bifurcation. Further, noting that the bifurcating branch emerges with
period
\begin{equation}
  \label{eq:omega}
  T=\tau\frac{2\pi}{\sqrt{\lambda^2S_1^2( b c - a d)+ \lambda S_1 (d-a)  + 1 }},
\end{equation}
and fixing $T=\frac{1}{2f}$ we can set
$\tau=\frac{\sqrt{\lambda^2S_1^2( b c - a d)+ \lambda S_1 (d-a) + 1 }}{4f\pi}$
such that the frequencies of oscillations produced at the Hopf match
the forcing frequency. 

\begin{figure}[t]
	\centering
{\includegraphics[width=120mm]{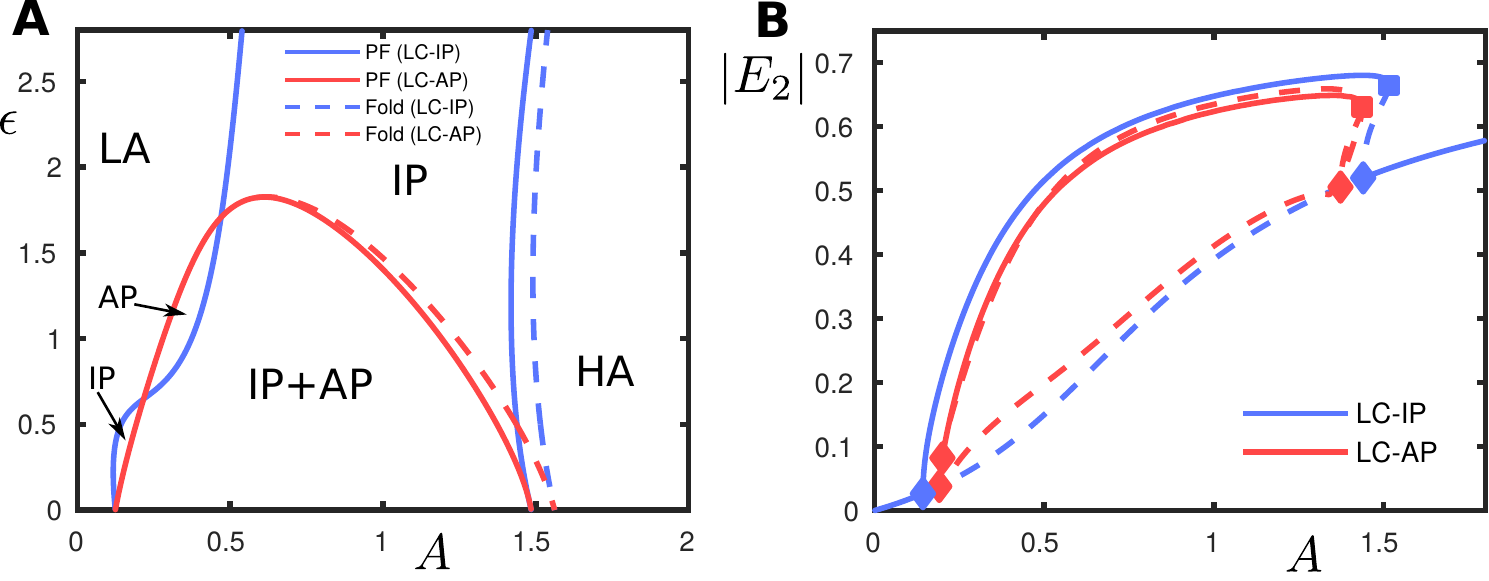}}
\caption{Bifurcation diagram with parameters $\mathcal{P}$ whilst
  setting $\lambda$ and $\tau$ as described in the text. \textbf{A}:
  Two-parameter bifurcation diagram in the $(A,\epsilon)$ plane
  showing locus of bifurcations with legends and labelling as in
  \ref{fig:bifFirstCase}; LA is a symmetric ($(E_1,I_1)=(E_2,I_2)$)
  low-amplitude limit cycle oscillations (following the periodic
  input) and HA is symmetric high-amplitude limit cycle. The IP and AP
  solutions co-exist in the region up to the dashed fold curve to the
  right. \textbf{B}: One-parameter bifurcation diagram at fixed
  $\epsilon=0.5$; dashed curve segments are unstable. Diamonds are
  pitchfork bifurcations and squares are fold bifurcations. The stable
  IP branch exists between a pitchfork bifurcation to the left and
  fold to the right. The AP branch emerges unstable and is stable
  between a secondary pitchfork bifurcation on the left and a fold
  bifurcation to the right.} \label{fig:bifPF}
\end{figure}

\Fref{fig:bifPF} shows a bifurcation diagram for the pair periodically-forced Wilson-Cowan oscillators. 
Each E-I oscillator receives the same input ($h=0$). Panel \textbf{A} shows regions of the $(\epsilon,A)$ plane in which different types of oscillatory behaviours are stable. For low forcing amplitude there are only low-amplitude oscillations, effectively modulating the FP solution in the unforced system. As $A$ is increased, Pitchfork bifurcations give rise to stable IP and AP branches that coexist (see panel B) for small $\epsilon$ approaching the uncoupling limit. For large $\epsilon$ the IP solution persists at intermediate values of $A$. For large $A$ there is a saturated high-amplitude solution.

The key result here is that the behaviour found in the unforced system is preserved for sufficiently small coupling strength and for weak forcing (IP and AP solutions persist close to the uncoupling limit, IP+AP region in \Fref{fig:bifPF}A). For larger forcing amplitude the intrinsic dynamics is overwhelmed and the forcing modulates a symmetrical fixed point (HA region in \Fref{fig:bifPF}A). This bifurcation analysis demonstrates the possibility for coexisting in-phase and anti-phase responses of the coupled Wilson-Cowan oscillators to encode network states corresponding to ``percept 1'' (IP) and ``percept 2'' (AP) as described in Section \ref{sec:intro}. This is possible without strong mutual inhibition (i.e. in the uncoupling limit) between abstract representations of the two possible percepts.

\section{Discussion and conclusions}

The study of identical coupled oscillators near a Hopf bifurcation is applicable to a wide range of systems where near-identical units undergo oscillatory instability. 
These systems may in general be represented by very different vector fields. 
Using the normal form theory of in \cite{ashwin2016hopf}, we are able to predict universal aspects of the mathematical behaviour for such systems. 
The analysis performed in this work for two oscillators reveals that, as is often the case in normal forms, 
although \eqref{eq:tereSystem} involves a big number of parameters, in the weak coupling limit, just a few of them govern and determine the possible bifurcations of the system.

Because of the symmetries of the system, there are usually two phase-locked oscillating solutions corresponding to in-phase ($\Delta \varphi = 0$) and anti-phase ($\Delta \varphi = \pi$). Depending on parameters, we find that all possible combinations between different stabilities of both solutions are possible. Our numerical analysis has shown that away from the coupling limit, richer dynamical behaviour is possible, with secondary bifurcations from the anti-phase branch and regions of coexistence between fixed-point and anti-phase solutions mediated by a fold of cycles. These scenarios can include modulated states that appear at torus bifurcations (see for example Figure~\ref{fig:bifSecCase}). Furthermore, we find the coexistence of in-phase and anti-phase solutions persists even in the presence of periodic forcing.

\subsection{Implications for models of perceptual bistability and neural competition}

Models of perceptual bistability are widely based on the assumption of strong mutual inhibition between populations of neurons that encode different perceptual interpretations of ambiguous stimuli. In general, this assumes that populations associated with different percepts are separated in some feature space (e.g. orientation in binocular rivalry) and that these populations enter into competition through mutual inhibition. However, when stimuli are periodic and the two possible perceptual interpretations involve the same features, it is less clear how competition between percepts might arise. For example, for the visual (auditory) stimulus in \fref{fig:schematic} both ``percept 1'' and ``percept 2'' involve the left spatial location (higher pitch A tone). It is therefore unclear how mutual inhibition between ``percept 1'' and ``percept 2'' could be implemented in neural hardware (although see \cite{rankin2015neuromechanistic} where a population pooling inputs from an intermediate feature location was proposed). Another possibility, proposed and demonstrated to be feasible in this study, involves oscillatory neural activity. Indeed, encoding of perceptual interpretations through oscillations allows for complete synchronisation of the network with all incoming inputs (like ``percept 1'') or for partial synchronisation of different parts of a network with separate elements (here in anti-phase). Furthermore, such an encoding mechanism does not rely on strong mutual inhibition, widely assumed between  the abstracted percept-based neural populations in competition models with little supporting evidence.

\subsection{Future perspectives}

An obvious extension of the bifurcation analysis would be to the forced symmetry broken case. If there is no assumed symmetry between percepts 1 and 2, this will result in a separation of Hopf bifurcations in the uncoupled limit and presumably mode locking and torus breakup scenarios familiar from the non-symmetric Hopf-Hopf interaction case \cite{gavrilov1980some}. Finally, one can consider the periodically forced system. Periodic forcing of the oscillators considered here (e.g. \cite{kim2015signal} for a single oscillator) will bring us to potentially much more complex bifurcation problems.

The study has demonstrated the potential role of oscillations in encoding different interpretations of periodically modulated ambiguous stimuli. It remains to explore the further role of feature space (say spatial location or tone frequency) and its interaction with oscillatory mechanisms. Additionally, as bistable perception involves spontaneous switching between perceptual interpretations, the mechanisms for these switches in the light of oscillatory stimuli remains to be explored.

Perceptual bistability with periodically modulated stimuli is robust over a range of input rates for the stimulus, whereas the simple network motif studied here has a fixed preferred input rate. So-called gradient networks of coupled oscillators have been proposed as a framework to understand many elements of early auditory processing and for perception of musical rhythm and beat \cite{large2010canonical,large2015neural}. Such a framework could be extensible to the study of perceptual bistability, relying in the dynamic mechanisms proposed here in the simple case of only two coupled oscillators.

\section*{Declarations}

\subsection*{Ethics approval and consent to participate}

Not applicable

\subsection*{Consent for publication}

Not applicable

\subsection*{Availability of data and material}

Not applicable

\subsection*{Acknowledgements}

This work has been partially funded by the grants MINECO-FEDER MTM2015-65715-P, MDM-2014-0445, PGC2018-098676-B-100 AEI/FEDER/UE, the Catalan grant 2017SGR1049, (GH, AP, TS), the MINECO-FEDER-UE MTM-2015-71509-C2-2-R (GH), 
and the Russian Scientific Foundation Grant 14-41-00044 (TS).
GH acknowledges the RyC project RYC-2014-15866. TS is supported by the Catalan Institution for research and advanced studies via an ICREA academia price 2018.
AP acknowledges the FPI Grant from project MINECO-FEDER-UE MTM2012-31714. We thank T. L\'azaro for providing us valuable references to compute the normal form coefficients.
We also acknowledge the use of the UPC Dynamical Systems group's cluster for research computing.\footnote{See: \href{https://dynamicalsystems.upc.edu/en/computing}{https://dynamicalsystems.upc.edu/en/computing}.} PA and JR acknowledge the financial support of the EPSRC Centre for Predictive Modelling in Healthcare, via grant EP/N014391/1. JR acknowledges support from an EPSRC New Investigator Award (EP/R03124X/1). 

\subsection*{Competing Interests}

The authors declare they have no competing interests.

\subsection*{Authors' Contributions}

PA, JR, AP formulated the problem. All authors were involved in the theoretical analysis, discussion of the results and writing the manuscript. AP and JR performed numerical simulations. All authors have read and approved the final version.

\bibliography{bibSym}
\bibliographystyle{abbrv}

\newpage
\appendix

\section{Appendix: Coupling terms}\label{sec:appendix}

The coupling terms $f_r$ and $f_\varphi$ for system \eqref{eq:polarEqs} are given by
\begin{equation}\label{eq:basicEqs}
\begin{split}
f_r (r_1, r_2, \Delta \varphi) & = r^2_1r_2 \left[(\beta_{\epsilon 1R} + \alpha_{\epsilon 3R})\cos(\Delta \varphi) - (\beta_{\epsilon 1I} - 
\alpha_{\epsilon 3I})\sin(\Delta \varphi)\right] \\ 
& + r^2_2r_1\left[\alpha_{\epsilon2R} + \beta_{\epsilon3R}\cos(2\Delta \varphi) - \beta_{\epsilon3I}\sin(2 \Delta \varphi)\right]  + r_1 \alpha_{\epsilon0R} + r^3_1 \alpha_{\epsilon1R}\\ 
& +  r^3_2\left[\beta_{\epsilon2R}\cos(\Delta \varphi) - \beta_{\epsilon2I}\sin(\Delta \varphi)\right] 
+ r_2\left[\beta_{\epsilon0R}\cos(\Delta \varphi) - \beta_{\epsilon0I} \sin(\Delta \varphi)\right],\\
\\
f_\varphi(r_1, r_2, \Delta \varphi) &= r^2_1r_2\left[ (\beta_{\epsilon1I} + \alpha_{\epsilon3I})\cos(\Delta \varphi) 
+ (\beta_{\epsilon1R} - \alpha_{\epsilon3R})\sin(\Delta \varphi)\right] \\ 
&+ r^2_2r_1\left[ \alpha_{\epsilon2I} + \beta_{\epsilon3I}\cos(2\Delta \varphi) + \beta_{\epsilon3R}\sin(2 \Delta \varphi)\right] + r_1 \alpha_{\epsilon0I} + r^3_1 \alpha_{\epsilon1I} \\ 
& + r^3_2\left[ \beta_{\epsilon2I}\cos(\Delta \varphi) + \beta_{\epsilon2R}\sin(\Delta \varphi)\right] 
+ r_2\left[ \beta_{\epsilon0I}\cos(\Delta \varphi) + \beta_{\epsilon0R} \sin(\Delta \varphi)\right].
\end{split}
\end{equation}

The coupling term $f_{\Delta \varphi}$ for system \eqref{eq:phaseDifEqs} is given by
\begin{equation}\label{eq:phaseDifTerms}
\begin{aligned}
f_{\Delta \varphi}(r_1, r_2, \Delta \varphi) &= f_{\varphi}(r_2, r_1, -\Delta \varphi)/r_2 - f_{\varphi}(r_1, r_2, \Delta \varphi)/r_1 = \\ &=
\left(r_1^2 - r^2_2\right) \left[ \alpha_{\epsilon2I} - \alpha_{\epsilon1I} + \beta_{\epsilon3I}\cos(2 \Delta \varphi)\right] - 2r_1r_2(\beta_{\epsilon1R} - \alpha_{\epsilon3R})\sin(\Delta \varphi) \\ &- \left(r_1^2 + r^2_2 \right) \beta_{\epsilon3R}\sin(2 \Delta \varphi) + \left( \frac{r^3_1}{r_2} - \frac{r^3_2}{r_1} \right) \beta_{\epsilon2I}\cos(\Delta \varphi) - \left( \frac{r^3_1}{r_2} + \frac{r^3_2}{r_1} \right) \beta_{\epsilon2R}\sin(\Delta \varphi) \\ & + \left( \frac{r_1}{r_2} - \frac{r_2}{r_1} \right)\beta_{\epsilon0I}\cos(\Delta \varphi) - \left(\frac{r_1}{r_2}
+ \frac{r_2}{r_1} \right)\beta_{\epsilon0R}\sin(\Delta \varphi).
\end{aligned}
\end{equation}

The coupling terms $g_s$, $g_d$ and $g_{\Delta \varphi}$ for system \eqref{eq:sdEqs} are given by
\begin{equation}\label{eq:sdFunctions}
\begin{split}
g_s(s, d, \Delta \varphi) &= f_r \left(\frac{s+d}{2}, \frac{s-d}{2}, \Delta \varphi \right) + f_r \left(\frac{s-d}{2}, \frac{s+d}{2}, -\Delta \varphi \right) = \\ &=
s(\cos(\Delta \varphi)\beta_{\epsilon0R} + \alpha_{\epsilon0R}) + d\sin(\Delta \varphi)\beta_{\epsilon0I} + \frac{s}{4}\left(s^2 + 3d^2 \right) (\beta_{\epsilon2R}\cos(\Delta \varphi) + \alpha_{\epsilon 1R}) \\ & + \frac{s}{4}\left(s^2 - d^2\right)\left[ (\beta_{\epsilon1R} + \alpha_{\epsilon3R})\cos(\Delta \varphi) + \alpha_{\epsilon2R} + \beta_{\epsilon3R}\cos(2\Delta \varphi)\right] \\ & + \frac{d}{4}\left(s^2 - d^2 \right)\left[ \beta_{\epsilon3I}\sin(2\Delta \varphi) - (\beta_{\epsilon1I} - \alpha_{\epsilon3I})\sin(\Delta \varphi)\right] + \frac{d}{4} \left(3s^2 + d^2 \right)\beta_{\epsilon2I}\sin(\Delta \varphi), \\
g_d(s, d, \Delta \varphi) &= f_r \left(\frac{s+d}{2}, \frac{s-d}{2}, \Delta \varphi \right) - f_r \left(\frac{s-d}{2}, \frac{s+d}{2}, -\Delta \varphi \right) = \\ &=
-d(\cos(\Delta \varphi)\beta_{\epsilon0R} - \alpha_{\epsilon0R}) - s\sin(\Delta \varphi)\beta_{\epsilon0I} - \frac{d}{4}\left(d^2 + 3s^2 \right) (\beta_{\epsilon2R}\cos(\Delta \varphi) - \alpha_{\epsilon 1R}) \\ & + \frac{d}{4}\left(s^2 - d^2\right)\left[ (\beta_{\epsilon1R} + \alpha_{\epsilon3R})\cos(\Delta \varphi) - \alpha_{\epsilon2R} - \beta_{\epsilon3R}\cos(2\Delta \varphi)\right] \\ & - \frac{s}{4}\left(s^2 - d^2 \right)\left[ \beta_{\epsilon3I}\sin(2\Delta \varphi) + (\beta_{\epsilon1I} - \alpha_{\epsilon3I})\sin(\Delta \varphi)\right] - \frac{s}{4} \left(3d^2 + s^2 \right)\beta_{\epsilon2I}\sin(\Delta \varphi), \\
g_{\Delta \varphi}(s, d, \Delta \varphi) &= f_{\Delta \varphi} \left(\frac{s+d}{2}, \frac{s-d}{2}, \Delta \varphi \right) = \\ &=
\beta_{\epsilon0I}\cos(\Delta \phi) \left(\frac{4sd}{s^2 - d^2}\right)- 2\beta_{\epsilon0R}\sin(\Delta \phi)\left(\frac{s^2 + d^2}{s^2 - d^2}\right) \\ & - \beta_{\epsilon2R}\sin(\Delta \phi)\left(\frac{(s^2 + d^2)^2}{(s^2 - d^2)} -  \frac{(s^2 - d^2)}{2}\right) + \beta_{\epsilon2I}\cos(\Delta\phi)\frac{2sd(s^2 + d^2)}{(s^2 - d^2)} \\ & - \beta_{\epsilon3R}\sin(2\Delta\phi)\frac{(s^2 + d^2)}{2}  - (\beta_{\epsilon1R} - \alpha_{\epsilon3R})\sin(\Delta\phi)\frac{(s^2 - d^2)}{2} \\ &+ (\alpha_{\epsilon2I} + \beta_{\epsilon3I}\cos(2\Delta\phi) - \alpha_{\epsilon 1I})sd.
\end{split}
\end{equation}

The terms for the Jacobian matrix in \eqref{eq:jacobianMatrixSD} are given by
\begin{equation}\label{eq:jacobianTermsCs}
\begin{aligned}
c^s_s &= \lambda + \epsilon(\alpha_{\epsilon0R} \pm \beta_{\epsilon0R}) + \frac{3s^2}{4} \Big(\alpha_{01R} + \epsilon (\alpha_{\epsilon1R} \pm (\beta_{\epsilon2R} + \beta_{\epsilon1R} + \alpha_{\epsilon3R}) + \alpha_{\epsilon2R} + \beta_{\epsilon3R})\Big),\\
c^d_d &= \lambda + \epsilon(\alpha_{\epsilon0R} \mp \beta_{\epsilon0R}) + \frac{s^2}{4} \Big(3\alpha_{01R} + 
\epsilon(3(\alpha_{\epsilon1R} \mp \beta_{\epsilon2R}) \pm (\beta_{\epsilon1R} + \alpha_{\epsilon3R}) - \alpha_{\epsilon2R} - \beta_{\epsilon3R})\Big), \\ 
c^d_{\Delta\varphi} &= \epsilon \left(-\frac{s^3}{4}(2\beta_{\epsilon3I} \pm (\beta_{\epsilon1I} - \alpha_{\epsilon3I}) \pm \beta_{\epsilon2I}) \mp \beta_{\epsilon0I}s\right), \\
c^{\Delta\varphi}_d &= -\alpha_{01I}s + \epsilon\left(s(\alpha_{\epsilon2I} - \alpha_{\epsilon1I} + \beta_{\epsilon3I} \pm 2\beta_{\epsilon2I}) \pm 4\frac{\beta_{\epsilon0I}}{s}\right),\\
c^{\Delta\varphi}_{\Delta\varphi} &= \epsilon \left(\frac{s^2}{2} (\mp (\beta_{\epsilon1R} - \alpha_{\epsilon3R}) - 2\beta_{\epsilon3R} \mp \beta_{\epsilon2R}) \mp 2\beta_{\epsilon0R}\right),
\end{aligned}
\end{equation}
where $s = s^{\pm}_{osc}$ and consistently with the notation used throughout the article, the $\pm$ sign corresponds to $\Delta \varphi = 0, \pi$ respectively.

\section{Normal Form Computation}\label{sec:normalFormCoefs}

In this appendix we provide a brief description of the numerical procedure used to compute the coefficients of the normal form \eqref{eq:tereSystem}. The procedure is related to normal form techniques in which one constructs a change of variables of the form
\begin{equation}\label{eq:normChangeCoord}
z = y + Q_2(y, \bar{y}) + Q_3(y, \bar{y}),
\end{equation}
where $Q_2(y, \bar{y})$ and $Q_3(y, \bar{y})$ are polynomials or order 2 and 3, respectively, such that the system \eqref{eq:originalSystem} expressed in the variables $y$ and $\bar{y}$ has the simplest expression possible. That is,
\begin{equation}
\label{eq:normChangeCoordDyn}
\dot{y} = Ay + f_2(y, \bar{y}) + f_3(y, \bar{y}) + \mathcal{O}_4(y, \bar{y}),
\end{equation}
where $A$ is the linearized sytem \eqref{eq:tereSystem} around the origin, $f_2 = 0$ and $f_3$ has the same monomials appearing in \eqref{eq:tereSystem}, namely: $y_1|y_1|^2$, $y_1|y_2|^2$, $y_4y_1^2$, $y_2|y_1|^2$, $y_2|y_2|^2$ and $y_3y_2^2$. 

To that aim we perform the following steps:

\begin{enumerate}
\item Consider the Taylor expansion of system \eqref{eq:originalSystem} around the origin:
\begin{equation}\label{eq:normalTaylor}
\begin{aligned}
\dot{z} &= Az + P_2(z, \bar{z}) + P_3(z, \bar{z}) + \mathcal{O}_4(z, \bar{z}),\\
\dot{\bar{z}} &= \bar{A}\bar{z} + \bar{P_2}(z, \bar{z}) + \bar{P_3}(z, \bar{z}) + \mathcal{O}_4(z, \bar{z}),
\end{aligned}
\end{equation}
where $z = (z_1, z_2)$, $\bar{z} = (\bar{z}_1, \bar{z}_2)$ $\in$ $\mathbb{C}^2$, $A = diag(\mu^+, \mu^-)$ is a diagonal matrix with $\mu^+, \mu^- \in \mathbb{C}$. $P_2$ and $P_3$ in \eqref{eq:normalTaylor} correspond to polynomials of degree 2 and 3, respectively. As $\bar{z}$ is the complex conjugate of $z$, we will just consider the first equation in \eqref{eq:normalTaylor}.

\item Compute the $q_{ij}$ coefficients of the polynomial $Q_2(y, \bar{y})$ given by
\begin{equation}\label{eq:qPolinomials}
Q_2(y, \bar{y}) = \Pi^{N=4}_{i=1}\Pi^{N=4}_{j=i}q_{ij}y_iy_j,
\end{equation}
where $y_3 = \bar{y}_1$ and $y_4 = \bar{y}_2$, by solving the following equation for each monomial
\begin{equation}
AQ_2(y, \bar{y}) - D_yQ_2(y, \bar{y}) Ay - D_{\bar{y}}Q_2(y, \bar{y}) \bar{A}\bar{y}
= f_2(y, \bar{y}) - P_2(y, \bar{y}). 
\end{equation}
With this choice, all the monomials in $f_2$ in \eqref{eq:normChangeCoordDyn} are null. 

\item Compute $f_3(y, \bar{y})$ given by the expression 
\begin{equation}
f_3(y, \bar{y}) = D_yP_2(y, \bar{y})Q_2(y, \bar{y}) + D_{\bar{y}}P_2(y, \bar{y})\bar{Q}_2(y, \bar{y}) + P_3(y, \bar{y}),
\end{equation}
thus obtaining the coefficients corresponding to the surviving monomials in \eqref{eq:tereSystem}: $y_i|y_i|^2$, $y^2_i\bar{y}_j$, $y_j|y_i|^2$, $y_i|y_j|^2$, $y^2_j\bar{y}_i$, $y_j|y_j|^2 \quad (i = 1,2 \quad j \neq i)$. 

\item Perform the change of coordinates $y = Cx$ in system \eqref{eq:normChangeCoordDyn}, where 
\begin{equation}\label{eq:matrixC}
C = \frac{1}{\sqrt{2}}
\begin{pmatrix}
1 & 1 & 0 & 0 \\
1 & -1 & 0 & 0 \\
0 & 0 & 1 & 1 \\
0 & 0 & 1 & -1 
\end{pmatrix},
\end{equation}
so that the system is written in the form \eqref{eq:tereSystem}.
\end{enumerate}

Notice that to compute the coefficients of $f_3$ in \eqref{eq:normChangeCoordDyn} it is enough to compute the change in \eqref{eq:normChangeCoord} up to order two. As a final remark, notice that apart from $\omega$ and $\alpha_{01}$ all the coefficients $\alpha_{\epsilon i},\beta_{\epsilon i}$ ($i = 0,...,3$) in \eqref{eq:tereSystem} are multiplied by  $\epsilon$. Therefore, to obtain the value of the coefficients we follow the procedure described above for $\epsilon = 0$, thus obtaining $\omega$ abd $\alpha_{01}$, and then repeat the same procedure for a small $\epsilon \neq 0$, which, using that $\omega$ and $\alpha_{01}$ are known, provides the coefficients $\alpha_{\epsilon i}, \beta_{\epsilon i}$. 

\end{document}